# GRACEFUL AND PRIME LABELINGS
## -Algorithms, Embeddings and Conjectures

**Suryaprakash Nagoji Rao***

Exploration Business Group
MRBC, ONGC, Mumbai-400 022

**Abstract.** Four algorithms giving rise to graceful graphs from a known (non-)graceful graph are described. Some necessary conditions for a graph to be highly graceful and critical are given. Finally some conjectures are made on graceful, critical and highly graceful graphs. The Ringel-Rosa-Kotzig Conjecture is generalized to highly graceful graphs. Mayeda-Seshu Tree Generation Algorithm is modified to generate all possible graceful labelings of trees of order p. An alternative algorithm in terms of integers modulo p is described which includes all possible graceful labelings of trees of order p and some interesting properties are observed. Optimal and graceful graph embeddings (not necessarily connected) are given. Alternative proofs for embedding a graph into a graceful graph as a subgraph and as an induced subgraph are included. An algorithm to obtain an optimal graceful embedding is described. A necessary condition for a graph to be supergraceful is given. As a consequence some classes of non-supergraceful graphs are obtained. Embedding problems of a graph into a supergraceful graph are studied. A catalogue of super graceful graphs with at most five nodes is appended. Optimal graceful and supergraceful embeddings of a graph are given.

Graph theoretical properties of prime and superprime graphs are listed. Good upper bound for minimum number of edges in a nonprime graph is given and some conjectures are proposed which in particular includes Entringer's prime tree conjecture. A conjecture for regular prime graphs is also proposed.

**Subjects**: Advanced Mathematics; Discrete Mathematics; Combinatorics;Graph Theory
**Keywords**: graceful and prime labelings; graceful graphs; Euler graphs;
**AMS subject classifications**: 05C78

## 1. INTRODUCTION

The word *'graph'* will mean a finite, undirected graph without loops and multiple edges. Unless otherwise stated a graph is connected. For terminology and notation not defined here we refer to *Harary (1972), Mayeda (1972), Buckly (1987)*.

A graph G is called a *'labeled graph'* when each node u is assigned a label $\varphi(u)$ and each edge uv is assigned the label $\varphi(uv)=|\varphi(u)-\varphi(v)|$. In this case $\varphi$ is called a *'labeling'* of G. Define $N(\varphi)=\{n \in \{0,1,..., q_0\}: \varphi(u)=n$, for some $u \in V\}$, $E(\varphi) = \{e \in \{1,2,..., q_0\}: |\varphi(u)-\varphi(v)|=e$, for some edge $uv \in E(G)\}$. Elements of $N(\varphi)$, $(E(\varphi))$ are called *'node (edge) labels'* of G with respect to $\varphi$. Their set complements are $\overline{N}(\varphi) = \{0,1,...,q_0\}-N(\varphi)$ and $\overline{E}(\varphi) = \{1,2,...,q_0\}-E(\varphi)$, where $q_0 = \max\{\varphi(u): u \in V\}$ with ' $^-$ ' is set complementation. Elements of $\overline{N}$ and $\overline{E}$ are called the *missing node and missing edge labels* with respect to $\varphi$. For any $e \in \overline{E}(\varphi)$ define, $R(e) = \{uv \notin E(G)$ and $|\varphi(u)-\varphi(v)|=e\}$. R(e) is the set of non-edges in G with edge labels e induced by $\varphi$.

A (p,q)-graph G is *'gracefully labeled'* if there is a labeling $\varphi$ of G such that $N(\varphi) \subseteq \{0,1,...,q\}$ and $E(\varphi) = \{1,2,...,q\}$. Such a labeling is called a *'graceful labeling'* of G. A *'graceful graph'* can be gracefully labeled, otherwise it is a *'non-graceful graph'* (see *Rosa (1967), Golomb (1972)*

---

*Sheppard (1976) and Guy (1977) for chronology 1969-1977 and Bermond (1978).* Bigraceful graph is a bipartite graceful graph. If a graph G is non-graceful, then a labeling φ which gives distinct edge labels such that the maximum of the node labels is minimum is called an *'Optimal labeling'* of G and the graph is called *'optimally labeled (optimal) graph'*. Note that $0 \varepsilon N(\varphi)$. This maximum is denoted by opt(G)≥q with equality holding whenever G is graceful. For every optimal labeling θ of G the labeling $\bar{\theta}$ which assigns opt(G)-θ(u) node label to u is again an optimal labeling of G called *'complementary node labeling'* of G. $\bar{\theta}$ preserves edge labels as that of θ and θ(u)=0, d(u)=1 implies that $\bar{\theta}(v)=0$ where uv∈E. A node u (an edge e=uv) of G is called *i- attractive (i-repelling)* if there (there does not) exist a graceful labeling α of G such that α(u)=i (|α(u)-α(v)|=i).

In order to understand the structural properties of graceful graphs it would be useful to describe graphs whose sub-graphs are graceful. This leads to the following definitions: A graph G is *'highly graceful'* if every connected subgraph of it is graceful (see [4]). A non-graceful graph G is *'critically non-graceful (or critical)'* if every connected subgraph of it is graceful. A graph G is *`v- (e-) minimal'* with respect to a property if G-u (G-e) possess the property for each u ε V(G) (e ε E(G)). Note that highly graceful and critical graphs imply v- and e- minimal properties.

A (p,q)-graph is *'totally labeled'* if there is a labeling of G such that $N(\varphi) \cup E(\varphi) = \{1,2,...,p+q\}$ and such a labeling is called a *'total labeling'*. A *'super graceful (or totally labeled) graph'* admits a total labeling, otherwise it is called a *'non-supergraceful (or semi-totally labeled) graph'* (see *Acharya [2]*). If G is non-supergraceful then a labeling φ for which $N(\varphi) \cup E(\varphi) \subseteq \{1,2,...,p+q+q_0\}$ for some $q_0 > 0$ and $q_0$ minimum is called a *'semitotal labeling'*. Study of supergraceful graphs is interesting in its own right because a graceful graph can always be constructed from a supergraceful graph. By *'G is embedded into H'* we mean that there is a subgraph of H isomorphic to G. Clearly, G is embedded into itself. Let μ be a labeling of H such that $\mu_V=\varphi$, where $\mu_V$ is the *restriction of μ to V*. Then μ is called an *'extension'* of φ with respect to G. Let $T^i$ be the tree obtained recursively deleting all endnodes at i-th step (i=0,1,…) from $T^{i-1}$, where $T=T^0$. λ be the step such that $T^\lambda=K_1$ or $K_2$. Then T is called a *λ-lobster*. Each tree is a λ-lobster for some λ≥0 and is a caterpillar for λ=2. Caterpillars are graceful. *Bermond (1979)* conjectured that 3-lobsters are graceful. An endnode of a caterpillar C is called an *extreme node* if it is adjacent to an end node of $C^1$. An end node of a tree T is *central* if the caterpillar containing it meets the tree at a node of center of T, which is either $K_1$ or $K_2$.

Essence of this paper was reported in a series of reports (*Hebbare [30 to 34]*) under the *Research Project No. HCS/DST/409/76* while the author was at Mehta Research Institute, Allahabad (1980-81). Main effort is to establish or propose general properties of (non-)graceful graphs.

Varied *applications* of labeled graphs have been cited in the literature, for example, *communication Networks, X-Ray Crystallography, Radio Astronomy, Circuit Layout Design and Missile Guidance*. See *Bloom (1977), Bloom and Golomb (1977, 1978), Harary (1988)* and the references therein. Labeled directed graphs were studied and applied to *Algebraic Systems, Generalized Complete Mappings, Network Addressing Problems And N-Queen Problems* by *Bloom, Proul and Hsu; See Harary (1988)*.

## 2. KNOWN NON-GRACEFUL GRAPHS

From every labeling $\varphi$ of a (p,q)-graph and a given integer $n \geq 1$ the node set V can be partitioned into classes $V_i$ consisting of node labels $\equiv i \pmod{n}$. We refer to this as and *n-ary partition* of G with respect to $\varphi$. Probe structural properties of graphs whose degrees satisfy $\equiv 0 \pmod n$, for $n>0$ especially when n is prime. Consider an n-ary partition of a graph G (V,E) into $V_0, V_1, ..., V_{n-1}$ where $V = \cup V_i$, $i=0,1,...,n-1$ and $V_i = \{x: x \equiv i \pmod n\}$. Let $e_i = |E(<V_i>)|$, $i=0,1,...,n-1$, $e_{ij}=|E(<V_i,V_j>)|$, $i \neq j$, $i,j=0,1,...,n-1$. $q_k = \Sigma e_{ij}$, $|i-j|=k$, $k=0,1,...,n-1$; when $i \neq j$, $e_{ii}=e_i$, where $E(<V_i,V_j>)$ is the set of edges joining a node of $V_i$ to a node of $V_j$.

**Theorem 2.1**. A necessary condition for G to be graceful is that there exists an n-ary partition of V into $V_0, V_1, ..., V_{n-1}$ so that $q_i \leq [(q+(n-i)]/n$, $i=0,1,...,n-1$.

For given p and q different n-partitions can be considered for $1 \leq n \leq q$. Graphs not satisfying such conditions are not known except for n=2 and such construction is open. When n=2 we have:

**Theorem 2.2.** *(Golomb)* A necessary condition for a graph to be graceful is that there exists a binary partition of V into $V_1=O$ and $V_2=E$ so that $E(E,O)=[q+1]/2$, where O, E are the odd and even label sets.

Binary labeling helps in proving some of the complete and near complete graphs to be non-graceful. For example, it follows that $K_p$, p = 5,7,8,10,12,13,14,15,17,... fail to give binary partitioning and are non-graceful. In fact the following is true for $K_p$.

**Theorem 2.3**. *(Golomb)* $K_n$ is non-graceful for $n \geq 5$.

**Theorem 2.4.** *(Rosa, Golomb)* A necessary condition for an eulerian (p,q)-graph to be graceful is that $[(q+1)/2]$ is even.

This is due to *Rosa (1967), Golomb (1972)* and an eulerian graph of this type is called here a *'Rosa-Golomb graph'*, shortly *'RG- graph'*. This implies that Eulerian graphs with $q \equiv 1 \text{ or } 2 \pmod 4$ are non-graceful. This prompts to classify eulerian graphs into four classes. Denote by $E_i$- the class of eulerian (p,q)-graphs with $q \equiv i \pmod 4$. For a graph G denote the number of sub-graphs of the type $E_i$ by $e_i$ for i=0,1,2,3. Most of the known non-graceful graphs contain a subgraph isomorphic to an RG-graph. *Golomb (1972)* also gives some other classes of non-graceful graphs. An isolated example of non-graceful graph is two copies of $K_4$ with a node in common (Golomb (1972)). $K_n$ is eulerian for n odd and non-graceful for $n \geq 5$. Probe to find new ones is continuing and interesting, especially the ones having small q.

**Problem 2.5.** Investigate the class of near complete graphs for non-gracefulness and graphs with path (cycle) of $K_n$s or $C_n$s that is, each node of the path (cycle) is replaced by a $K_n$ or $C_n$, $n \geq 0$.

Consider the family of graphs $H(l,m,n) = mK_l \times \overline{K_n}$. It was conjectured *Hebbare* (1976) that $H(2,2,n)$, $n>1$ and n even are non-graceful. Note also that they are not completely of RG- type. *Bhat-Nayak and Gokhale (1986)* have established the conjecture affirmatively.

**Theorem 2.6.** (*Bhat-Nayak & S.K. Gokhale*) $H(2,2,n)$ is non-graceful, $n>1$.

**Problem 2.7.** Investigate the family of graphs $H(l,m,n)$ for l, m, n >1 for non-gracefulness.

## 3. HIGHLY GRACEFULNESS AND CRTICALITY

Definition of a highly graceful graph G requires that G itself is graceful. Clearly, any non-graceful graph of order p is a forbidden subgraph of a highly graceful graph of order $\geq p$. In particular, an eulerian (p,q)- graph G with $q \equiv 1 or 2 \pmod 4$ is non-graceful and hence forbidden for all graphs of order $\geq p$. It is well known that paths are graceful. So a cycle $C_n$, $n \equiv 0 or 3 \pmod 4$ is highly graceful. Whereas $C_n$, $n \equiv 1 or 2 \pmod 4$ are non-graceful and so are critical. Describing eulerian graphs that are highly graceful or critical with minimal property would be interesting.

**Proposition 3.1.** G is highly graceful then the following graphs are forbidden:
1. A cycle $C_n$, $n \equiv 1 or 2 \pmod 4$,
2. Two cycles $C_m$, $C_n$, $m,n \equiv 3 \pmod 4$ with exactly one node in common.
3. Three cycles $C_l$, $C_m$, $C_n$ such that $l \equiv 0 \pmod 4$, $m,n \equiv 3 \pmod 4$ and $C_l$ has exactly one common node with each of the cycles $C_m$ and $C_n$.

That $C_m$ and $C_n$ are disjoint in (3) follows from (2). The graphs satisfying 1), 2) and 3) in Proposition 3.1 are non-graceful. (Non-eulerian) Sub-graphs of interest in general are as follows:

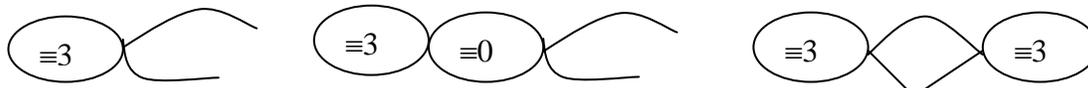

**Conjecture 3.2.** Graphs in 2) and 3) of Proposition 3.1 are critical.

The simplest class of connected graphs is trees. Even this class exhibits complexity with respect to labelings. Attempts to show them to be graceful lead to the following:

**Conjecture 3.3.** (*Ringel, Rosa, Kotzig*) All trees are graceful.

We call it RRK-conjecture and if true implies that all trees are highly graceful. Further, the truth of the conjecture implies that every non-graceful graph has a cycle. For a detailed exposition on this conjecture we refer to *Golomb (1972), Bloom (1979), Kotzig (1973), Rosa (1967), Bermond (1978)* and the references cited therein. That paths $P_n$ are graceful, for all $n \geq 1$ is well known. A strong and beautiful result is:

**Theorem *3.4*.** (*Flandrin, Fournier and Germa (1983)*) Each node of a path $P_n$, $n \geq 8$ is i-attractive for $i=1,\ldots,n$.

This was a conjecture of Cahit. It follows from this that cycles $C_n$, $n \equiv 0 or 3 \pmod 4$ are highly graceful whereas, cycles $C_n$, $n \equiv 1 or 2 \pmod 4$ are non-graceful and hence are critical. A natural question is '*Are there other graphs with property as in Theorem 3.4? If not describe i-attractive nodes for all $i=0,1,\ldots,q$ node labels*'. We propose the following:

**Conjecture 3.5.** Every extreme node of a caterpillar is i-attractive for $i=0,1,\ldots,q$.

Investigate Theorem 3.4 for edge labels of $P_n$, which leads to the following:

**Problem 3.6.** Is every edge of $P_n$ i-attractive for $i=1,\ldots,n$?

**Proposition 3.7**. A necessary condition for highly graceful is that RG-graphs are forbidden.

Non-graceful (p,q)- graphs exist for each p≥5. Non-graceful graphs of orders p≤6 are $C_5$, $K_5$; $C_6$, $K_6$ and $K_4$ with two nodes of degree two adjacent to disjoint pairs of nodes of $K_4$. v- (e-) minimal graphs with this property are of interest. For a non-graceful (p,q)- graph denote the minimum q by $\chi_p$. Within among eulerian graphs, cycles $C_n$, n≡1or2(mod 4) are non-graceful and have minimum number of edges. Such graphs for orders ≡0or3(mod 4) are again of interest. A (4t,4t+2)- graph t≥2 consisting of a $C_{4(t-1)}$ and two disjoint $K_3$s each having exactly one common node with $C_{4(t-1)}$ can be verified to be eulerian and so non-graceful. Further, a (4t+3, 4t+5)- graph t≥2 consisting of a $C_{4t-1}$ and two disjoint $K_3$s each having exactly one common node with $C_{4t-1}$ can be verified to be eulerian and so non-graceful. For p=7, we have the following whose proof will be given elsewhere.

**Proposition 3.8.** An eulerian graph of order 7 has $C_5$, $C_6$ or 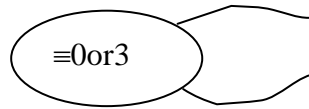 as a subgraph and so an eulerian critical graph of order 7 does not exist.

**Proposition 3.9.** For eulerian graphs: $\chi_p$=p, if p≡1or2(mod 4), p≥5. $\chi_p$ =p+2, if p≡0or3(mod 4), p≥8. When p=7 $\chi_p$ does not exists.

**Proposition 3.10.** $C_n$ is critical and so e-minimal for n≡1or2(mod 4).

Cycle structure of eulerian graphs is of interest. In particular, cycles $C_m$ and $C_n$ which intersect in a path $P_l$, l>0 give rise to a third cycle of length as in Table 1 with '≡' meaning ≡ (mod 4):

| Table 1. | l<br>m↓ n→ | ≡0 or 2 | | | | ≡1 or 3 | | | |
|---|---|---|---|---|---|---|---|---|---|
| | | ≡0 | 1 | 2 | 3 | ≡0 | 1 | 2 | 3 |
| | ≡0 | 0 | 1 | 2 | 3 | 2 | 3 | 0 | 1 |
| | ≡1 | 1 | 2 | 3 | 0 | 3 | 0 | 1 | 2 |
| | ≡2 | 2 | 3 | 0 | 1 | 0 | 1 | 2 | 3 |
| | ≡3 | 3 | 0 | 1 | 2 | 1 | 2 | 3 | 0 |

$E_{03}$ be the class of eulerian graphs with cycles of type ≡0or3(mod 4). In particular, consider a subclass $E'_{03}$ of $E_{03}$ such that two cycles $C_m$ and $C_n$ in G in $E'_{03}$ intersect exactly in a path $P_l$, l≥1implies that:
- l≡0or2(mod 4) and either m,n≡0(mod 4) or m≡0(mod 4) and n≡3(mod 4).
- l≡1or3(mod 4) and both m,n are ≡3(mod 4).

**Conjecture 3.11.** Graphs of $E'_{03}$ are highly graceful.            ≡0or3

Further, if cycles of G are all of type ≡0(mod 4) then from Table 1 it follows that any two intersecting cycles in $P_l$ satisfy l≡0or2(mod 4). Such a graph is bipartite. Additionally if G is highly graceful then cycles $C_n$, n≡2(mod 4) are forbidden and so a special case of Conjecture 3.11 is:

**Conjecture 3.12.** A graph G whose cycles $C_n$ are of length n≡0(mod 4) is highly bigraceful.

This generalizes RRK- conjecture. A large family of graphs with this property is known to be graceful. (See *Lee and Wang [41], Koh, Rosers and Lim (1979), Koh, Rosers, Lee and Toh [39], Koh, Rosers and Tan (1980)* and papers on graceful trees. Conjecture 3.10 was further generalized to d-sequential graphs in *Acharya (1983)*. Some classes of graphs of this type follow:

**Construction 3.13.** Let T be a tree with α-labeling and bipartition $V(T)=A\cup B$ with A having the labels $0,1,\ldots,|A|$. Construct the graph $T_B$ with $V(T_B)=A\cup B\cup B_1\cup\ldots\cup B_n$, $|B|=|B_1|=\ldots=|B_n|$ and $E(T_B)=E(T) \cup\{ab_i$ is an edge in G whenever ab is an edge in T for $a\in A$, $b\in B$, $b_i\in B_i\}$, $i=1,\ldots,n$. $T_B$ is union of T with n copies of partition B so that $<A\cup B_i> \sim T$. $T_B$ is bipartite and has only 4-cycles. Also if G is a graph with cycles $n\equiv 0\pmod 4$ then a 1-subdivision graph is again of that type. Truth of conjecture 3.11 implies that unicyclic graphs with the cycle $C_n$, $n\equiv 0\pmod 4$ are bigraceful. But, we make the following general conjecture for unicyclic graphs and support it by giving a unicyclic algorithm to construct a large class of unicyclic graphs.

**Conjecture 3.14.** The only non-graceful unicyclic graphs are cycles $C_n$, $n\equiv 1\text{ or }2\pmod 4$.

If proved true implies that the only non-graceful, non-bipartite unicyclic graphs are the cycles $C_n$, $n\equiv 1\pmod 4$ and is verified for unicyclic graphs of order $\leq 8$ and all unicyclic graphs of order 9 with a 5-cycle, see also *Truszczynski (1984)*. An interesting class to probe into is unicyclic graphs having caterpillars as components after deleting nodes of the cycle. This class can be further generalized to have components as α-valuable trees. That unicyclic graphs with a caterpillar at exactly one node of the cycle is graceful follows from complementary node labeling.

**Construction 3.15.** Construction 1 can be modified so that each $B_i$ is subset of B with the edges as in T. In particular, if n=1 and $B_1$ has only one node of degree 2 then $T_B$ is a unicyclic graph.

**Construction 3.16.** Let G be a (p,q)-graph and φ, φ' be two α-labelings (possibly φ=φ') of G with A, B as partition of V(G) and $\varphi(A)=\{0,1,\ldots,|A|\}=\varphi'(A)$. H be the graph obtained from G, G' labeled graphs by φ, φ' so that the nodes of A of G, G' are identified whenever the nodes have the same label. The node labels of B' are added q. That is, $V(H)=A\cup B\cup B'$, $E(H)=E(G) \cup\{E(G')\}$. H is a (p+|B|,2q)-graph and is graceful with an α-labeling.

**Construction 3.17.** *'Panel graph'* consists of $l\geq 2$ non-empty subsets of isolated nodes with a connected bipartite graph between consecutive subsets. This can be further generalized to a *'Fence graph'* having several panels with 2 or more panel graphs meeting at a subset of another panel graph. *'Complete panel' ('Complete Fence') graph* has complete bipartite graph in between two consecutive subsets. Complete Panel graphs (sequential join of subsets of isolated nodes) are shown to be graceful (*Lee and Wang [41]*). We feel that complete fence graphs and in general fence graphs are also graceful. One subclass of interest is fence with graph induced in between two consecutive subsets is a tree. Further, when this induced graph is a path the panel is a partial grid graph. A graceful labeling can easily be given for such panels starting with an α-labeling of a path (see also *Acharya [2]*).

**Conjecture 3.18.** An eulerian graph G for which RG- graphs are forbidden is highly graceful.

If proved true implies that a bipartite eulerian graph G without cycles $C_n$, $n\equiv 2\pmod 4$ is highly graceful. This essentially means that no other types of non-graceful bipartite eulerian graphs are possible. Unlike complete and near complete graphs are more likely to be non-graceful, complete bipartite graphs are known to be graceful and near complete bipartite graphs are likely to be graceful. $B_{n,n}=K_{n,n}$-{perfect matching} is a $(2n,n(n-1))$-graph and is a regular bipartite graph of degree n-1. When n is odd and $n(n-1)\equiv 1 or 2 \pmod 4$ $B_{n,n}$ is eulerian and is non-graceful. For example $B_{7,7}$ is non-graceful. $B_{n,n}$ is known to be graceful for $n \leq 9$ except n=8.

**Problem 3.19.** Investigate $B_{n,n}$, $n \geq 10$ for gracefulness.
For bipartite, critically non-graceful graphs we propose the following:

**Conjecture 3.20.** The only bipartite, critical graphs are the cycles $C_n$, $n \equiv 2 \pmod 4$.

We state correct version of the misprinted conjecture in *Gangopadhyay and Rao Hebbare (1980)* and propose stronger version. By '0% of graphs with a property in a class of graphs' we mean the ratio of number of graphs with property of order p to the total number of graphs in the class of order p tends to zero with p tending to infinity.

**Conjecture 3.21.** 0% of all bipartite non-eulerian graphs with $q \equiv 2 \pmod 4$ are non-graceful.

**Conjecture 3.22.** 0% of all bipartite graphs are non-graceful.

The concepts of highly graceful and critically non-graceful graphs aided in solidifying the ideas and manifestation of the structure of graceful graphs. Series of conjectures in the previous sections are the result of this analysis and can be represented schematically. Possible (expected) highly gracefulness and possible (expected) criticality are shown **Figs. 1a,b (2a,b)** with hashed regions indicating graphs without (with) the property. No meaning is assigned to region areas. In fact, the author's strong conviction is the following that extends Conjecture 3.18:

**Conjecture 3.23.** A graph G is highly graceful if and only if RG- graphs are forbidden in G.

## 4. GRACEFUL ALGORITHMS

An *algorithm* can be described as a system which accepts some graphs, functions, parameters, etc. as input and yields an output in terms of graphs, functions, parameters, etc. as illustrated in the diagram below. In other words, algorithms are analogous to filters in Digital Signal Processing. An algorithm may be implemented and automated in any computer language.

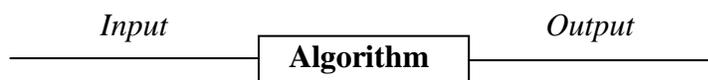

A *'graceful algorithm'* is one which when performed on a graceful or non-graceful graph as input gives rise to a graceful graph as output. Such algorithms are not new in the literature, for example see [1, 2, 4, 6, 7, 12, 13, and 14] and the references cited therein.

Four general graceful algorithms now follow. These algorithms are existential in the sense that new graceful graphs are constructed from a given graceful or non-graceful graph exhibiting an optimal labeling, subject to the existence of certain nonadjacent pairs of nodes in G or existence of a graph (not necessarily connected) with a suitable node labeling. Though algorithms given work for non-graceful graphs with appropriate modifications, we confine to only graceful graphs.

**Observation 4.1.** Let $\varphi$ be a graceful labeling of a (p,q)-graph G so that $\varphi(u)=0$. Then one can mount a caterpillar at u and each farther extreme node from u of the caterpillar is 0-attractive. Conversely, if u is 0-repelling in G then each of the extreme nodes of the caterpillar at u is 0-repelling in the new graph. This follows from complementary labeling.

Algorithms discussed here are not efficient as they scan the power set or its subsets and are of complexity of order $2^n$ where n is the cardinality of set under consideration. However, depending on the p, q of a (p,q)-graph the number of steps might be small. Computationally efficient algorithms are also discussed for special classes of graphs. See *Mayeda (1972) and Chen (1971).*

**Algorithm 4.2.**
*Input.* A (p,q)-graph G and $\varphi'$- a graceful labeling of G.
*Output.* Graceful (p+m,q')-graph H with a graceful labeling $\gamma$ where m>0 and q'>q.

*Step 1.* Find m positive integers $c_1,...,c_m$, (1<m<p) such that a new labeling $\varphi$ results as follows: $\varphi(u_i) = \varphi'(u_i) + c_i$, $i=1,...,m$, where $u_1,...,u_m$ are so that $N(\varphi)$ and $E(\varphi)$ are distinct and $q_0=q+c$ with $c=\max\{c_i, i=1,...,m\}$. (When this maximum is attained for some i, say $u=u_i$, then either u must be adjacent to the node labeled 0 if it exists in $\varphi$ or the new graph should be of order greater than p with a node labeled 0 adjacent to u).

*Step 2.* Let $S \subseteq \overline{E}(\varphi)$ and $T \subseteq \overline{N}(\varphi)$. Find a graph D with k ($\geq 1$) components $D_1,...,D_k$ and a labeling $\gamma'$ of $D_i$ such that $N(\gamma'_i) \cap N(\varphi') = \{v_i\}$, for some $v_i \in V$, i=1,...,k not necessarily distinct and $|N(\gamma'_i) \cap N(\gamma'_j)| \leq 1$, for $i \neq j$; $i,j=1,...,k$. Here $\gamma'_i$ stands for restriction of $\gamma'$ with respect to $D_i$ satisfying $N(\gamma') = T$ and $E(\gamma') = S$.

*Step 3.* Find G(S,T)- the set of all graphs D with a node labeling $\gamma'$ for each pair (S,T) as above. (An edge label $e \in \overline{E}(\varphi)$ is of *'type 1'* if a pair of nodes from R(e) is chosen and are joined by an edge in G. It is of *'type 2'* if for a choice of S and T there exists $D \in G(S,T)$ and a labeling $\gamma'$ of D such that $e \in E(\gamma')$). Clearly each $e \in \overline{E}(\varphi)$ is either of type 1 or of type 2. If R(e) is nonempty for an edge $e \in \overline{E}(\gamma)$ then construct a graph H by joining the pair of nonadjacent nodes selected. Then $\gamma=\varphi$ is a graceful labeling and H is a graceful graph. (For a choice of $D \in G(S,T)$ and a labeling $\gamma'$ of D, either $M= \overline{E}(\varphi)-E(\gamma')$ is empty or not.)

*Step 4.* If M is empty then construct a graph H by coalescing the graphs G and $D_i$ at $v_i$ (i=1,...,m). (Again the labeling $\gamma$ is a graceful labeling of H and is graceful.) If M is not empty then for each $e \in \overline{E}(\varphi)-E(\gamma')$ choose a pair of nodes for each such e from R(e) and join them in G. (Note that R(e) is empty only for e=q).

Denote by *H*(G, $\varphi$) the set of all graceful graphs obtained using the above algorithm. If each $e \in \overline{E}(\varphi)$ is of type 1 then H is a (p,q+c)-graph otherwise H is of order greater than p with q+c edges. In particular, G is an induced subgraph of H whenever each $e \in \overline{E}(\varphi)$ is of type 2. **Fig. 3** illustrates Algorithm 4.2. Three algorithms that are special cases now follow:

**Algorithm 4.3.**
*Input.* A bigraceful (p,q)-graph G with the bipartition $V=A \cup B$ and $\varphi'$ is a graceful labeling of G

such that $0 \in A$ and $\varphi'(u) < \varphi'(v)$ holds for any $u \in A$ and $v \in B$.
*Output.* A graceful graph $H(G, \varphi)$.

Consider a labeling $\varphi$ of G as follows: $\varphi(x) = \varphi'(x) + c$, for $x \in B$, and $\varphi(x) = \varphi'(x)$, for $x \in A$. (Then $\bar{N}(\varphi) = \{m \in \{0,1,...,c+q\}: m \notin E(\varphi')\}$ and $\bar{E}(\varphi) = \{1,2,...,c\}$). As in Algorithm 4.2 find a graph D (possibly disconnected) and a labeling $\gamma'$.

**Example 4.4.** Let $\varphi'$ be a graceful labeling of a caterpillar C using bipartition technique *(Rosa (1967))*. If $V = A \cup B$ is the bipartition of C and $\varphi'$ is such that $0 \in A$ then $N_{A(\varphi')} = \{0,1,...,m-1\}$ and $N_{B(\varphi')} = \{m, m+1,...,p-1\}$. Consider $\varphi$ as in Algorithm 4.3. Then $\bar{N}(\varphi) = \{c+1,...,p+c-1\}$ and $\bar{E}(\varphi) = \{1,...,c\}$. Then for any $c < |A|$ join the pairs $(0,1),...,(0,c)$ in G with respect to $\varphi$. Then $\varphi$ can be verified to be a graceful labeling and so G is graceful. Further let T be a tree of order $c+1$ with a graceful labeling. Define $\mu$ as follows: $\mu(x) = \mu'(x) + (c+1)$, $x \in V(T)$. Then $\mu$ has exactly one node label, say for the node $u \in V(C)$, in common with $\varphi$. Coalesce C and T at u and H be the resulting graph. The labeling $\mu$ such that $\gamma(V(C)) = \varphi$, and $\gamma(V(T)) = \mu$, is a graceful labeling of H.

**Example 4.5.** Algorithm 4.3 can be altered to give $\varphi$ considering a (p,q)-graph G as in the algorithm with a graceful labeling $\varphi'$ of G as follows: $\varphi(x) = \varphi'(x) + c$, for each $x \in A$ and $\varphi(x) = \varphi'(x)$, for each $x \in B$ where $c \geq q+1$ is any positive integer. In particular, if we choose G to be a caterpillar C of order p as in Example 1, then graceful graphs of order $> p+1$ can be constructed. When the order is $p+1$ then the graph H obtained by adding a new node u adjacent to each of the nodes in H with labeling such that $\gamma(u) = 0$ and $\gamma(V(C)) = \varphi'$ is graceful.

**Fig. 4** illustrates Algorithm 4.3. Note that if G is a tree, D is a forest and edges of the type 1 are not used then H is a tree. When H is a tree Algorithm 4.3 is called a *'Tree Algorithm'*.

**Algorithm 4.6.**
*Input.* T is a tree of order p and $\varphi'$ is a graceful labeling of T. Further, $u,v \in N(T)$ with $\varphi'(u) = 0$ and $\varphi'(v) = p-1$ (Then $uv \in E(T)$).
*Output.* A class of graceful graphs.

Define a labeling $\varphi$ of T so that $\varphi(x) = \varphi'(x)$, for $x \neq v$ and $\varphi(v) = \varphi'(v) + c$, where $c \leq p-2$. Find $H(T, \varphi)$ as in Algorithm 4.3. Then any $H \in H(T, \varphi)$ is a $(p, p+c-1)$-graph and is graceful. (Note that graphs constructed are of same order p). **Fig. 5** illustrates Algorithm 4.6.

**(Unicyclic) Algorithm 4.7.**
*Input.* T is a tree of order p and $\varphi'$ is a graceful labeling of T.
*Output.* Each H(u,v) is a graceful unicyclic graph of order p.

Consider $\varphi$ as in Algorithm 4.6 and set $c=1$. Clearly, the only unused node label in $\varphi$ is $m = (p-1) - \min\{\varphi(x): x \neq u \text{ and } vx \in E(T)\}$. Find $y,z \in V(T)$, so that $yz \notin E(T)$ and satisfying $|\varphi(y) - \varphi(z)| = m$. H(y,z) graph is obtained from T by adding the edge yz in T. **Fig. 6** illustrates Algorithm 4.7 for the possible pairs: $\{(0,4),(1,5),(2,6),(3,7),(5,9)\}$.

**Example 4.8.** Consider $\alpha$- labeling of a tree T with a bipartition $V(T) = A \cup B$ with $0,1,...,[n/2]$ labels of nodes in A and $[n/2]+1,...,n$ that of B. Add 1 to the labels of B. Then the edge labels $n, n-1,...,2$ appear except 1. This can be obtained by adding an edge between any two consecutive node labels. Add 1 to all labels of B except n then al edge labels except 3 can be obtained. This gives unicyclic graph with the cycle a 5-cycle. Similar operations by adding 1 to each label of A except 0 and adding 1 to each label of A except 0 and 1, and so on may be considered.

## 5. GRACEFUL EMBEDDINGS

That every graph can be embedded in a graceful graph is established in *Bloom (1975)* and *Acharya [2]* and hence there exists no forbidden subgraph characterization for graceful graphs.

**Embedding Algorithm 5.1.** Embed an arbitrary graph into a graceful graph.
*Input.* Arbitrary graph G.
*Output.* A graceful graph with G as a subgraph.
- Label nodes of G so no edge label is repeated with 0 assigned to some node.
- Assign edges with missing labels. If x is a missing two cases arise:
  - if x is a node label, draw the edge from 0 to x
  - else add a new isolated node adjacent to the node labeled 0 to the graph and label it x.

We now present a short proof in the first case and exhibit an embedding with smaller number of nodes in the later case.

**Theorem 5.2.** Every graph can be embedded into a graceful graph.

**Proof.** Let G be a (p,q)-graph.. If G is non-graceful consider the complete graph $K_p$ and an optimal labeling $\varphi$ of $K_p$. Clearly, $0 \in N(\varphi)$. Let $u \in V$ be such that $\varphi(u)=0$. Further, if $e \in \bar{E}(\varphi)$ the $e \notin N(\varphi)$. Hence, construct a graph H adjoining $|\bar{E}(\varphi)|$ new nodes adjacent to u, one for each $e \in \bar{E}(\varphi)$ with the node weight e. Clearly H is a graceful graph. Now since any graph G of order p is a subgraph of $K_p$, H is a graceful embedding of G and the result follows. ∎

The graph G in the above theorem need not necessarily be connected. However, H is a connected embedding of G. The above result can be generalized to k- sequential graphs.

**Theorem 5.3.** Every graph can be embedded into a graceful graph as an induced subgraph.

**Proof.** Let G be a (p,q)-graph. If G is non-graceful consider an optimal labeling $\varphi$ of G. Clearly, $0 \in N(\varphi)$. Let $u \in V$ be such that $\varphi(u)=0$. Define, $Q_1(\varphi) = \{e \in \bar{E}(\varphi): e \notin N(\varphi)\} = \{e_1, e_2, ..., e_l\}$, $m = \max\{e: e \in \bar{E}(\varphi) - Q_1(\varphi)\}$, $M(\varphi) = \max\{\varphi(x): x \in V\}$, and v be the node such that $\varphi(v) = \max\{\varphi(x): \varphi(x)+m > M(\varphi)\}$. Let $\bar{E}(\varphi) - Q_1(\varphi) = \{d_1, d_2, ..., d_k\}$ where l+k=opt(G)-q and $1 \leq d_1 \leq d_2 \leq ... \leq d_k = m$. Construct a (p+l+m, q+l+m+k)- graph H containing G as a subgraph with a labeling $\mu$ as follows: Let $u_1, u_2, ..., u_{l+m}$ be new nodes each of which is adjacent to u. $\mu$ is such that $\mu_V = \varphi$, $\mu(u_i) = e_i$, i=1,2,...,l and $\mu(u_{l+j}) = \varphi(v)+j$, j=1,2,...,m. For each $d_i$ let $\mu(u_{l+d_i}) = \varphi(v)+d_i$ and $u_{l+d_i}$ be adjacent to v. It can now be verified that $\mu$ is a graceful labeling and hence H is a graceful graph with G as an induced subgraph. ∎

We now present an algorithm to obtain an optimal graceful embedding H of a graph G with respect to a given optimal labeling of G such that G is an induced subgraph of H with minimum possible nodes.

**Algorithm 5.4.**
*Input.* Let G be a non-graceful (p,q)-graph with an optimal labeling φ.
*Output.* Optimal graceful embedding H of graph G.

Consider $\overline{N}(\varphi)$ and $\overline{E}(\varphi)$ which are nonempty. For any subset $S \subseteq N(\varphi)$ and $T \subseteq (\varphi)$ find a graph F (not necessarily connected) and V(F)∩V is induced in G with a labeling μ such that each label of $\overline{E}(\varphi)$ occurs exactly once where $\mu_S = \varphi_S$. Let G(S,T) be the set of all graphs thus obtained. Further let $G(\varphi) = \cup G(S,T)$, for $S \subseteq N(\varphi)$, $T \subseteq P(\varphi)$. Any graph $H \in G(\varphi)$ is an optimal embedding of G with respect to φ. Further, H is of order at least p+1 and has opt(G) edges with G as an induced subgraph. Should H as above with minimum order be connected? This is shown to be true under some conditions.

**Observation 5.5.** If H is disconnected then let C be a component of it containing G. Then C contains a node of G labeled zero. Let $C_1$ be any other component of G. Subtract the minimum node label of $C_1$ from each of its node labels. Clearly 0 occurs as node label in this new labeling of $C_1$. Construct a new graph merging the component $C_1$ with the component C at those nodes, which have the same node labels. Repeat this to each of the components of H different from C. We get a connected graceful embedding with minimal possible nodes, a contradiction. In general, G need not an induced subgraph of H. Suppose for each component $C_1 \neq C$ of H if there is no pair of nonadjacent nodes in G but nonadjacent in $C_1$ with same node labeling then G remains an induced subgraph of a connected embedding with minimum nodes, a contradiction.

Define for each such optimal labeling φ of G *'optimal order'* as follows: opt(p) = min{|H|: $H \in G(\varphi)$, φ is an optimal labeling of G}. Any embedding $H \in G(\varphi)$ of order opt(p) is an optimal graceful embedding of G. We end this section with:

**Conjecture 5.6.** Every embedding $H \in G(\varphi)$ of order opt(p) is connected.

**Problem 5.7.** $K_2$ is graceful with θ(u)=0, θ(v)=1 and the labeled graphs with θ and θ' are isomorphic. Describe the graphs with this property that is G is label isomorphic with respect to θ and θ'. Is this the only graph with this property?

## 6. TOWARDS GRACEFULNESS OF TREES!

**Algorithm *6.1.*** *(Modified Mayeda-Seshu Algorithm).* Here we address to the problem of generating all graceful labelings of all trees of order n+1. Consider $K_{n+1}$ with its nodes labeled by 0,1,...,n. The edge labels {1,2,...,n} arise out of the edges: 0n,0(n-1),...,01,1n,1(n-1),...,12,...,(n-1)n, where edge labeled i (1≤ i ≤ n) is repeated (n-i+1)- times as 0i, 1(i+1),..., (n-i)n. Note that a graceful graph of order n including trees can be embedded in $K_n$. Let the reference tree in the MS- algorithm be chosen as follows: $t_0$ be $K_{1,n}$ with the n degree node labeled 0 and arranging the nodes cyclically with labels 1,...,n. following the symbolism in *Mayeda (1972)* let $t_1$, $t_2$ be two trees of $K_{n+1}$ at a distance 1. Let $t_1-t_2=(ij)$ and $t_2-t_1=(kl)$. Then $t_2=t_1 \oplus (ij,kl)$. We stipulate the condition that |ij|=|kl| then it follows that if $t_1$ is graceful then $t_2$ (need not be a tree) is also graceful. Further, if removal of edge ij and inclusion of edge kl results in a tree then $t_2$ is a graceful tree. This is assured if incoming edge is in $S_e(t_1)$- fundamental cut-set.

$T^{e1,...,ek} = \{t' \oplus (e_k b) : b \in Se_k(t_0) \cap Se_k(t_k'), t' \in T^{e1,...,ek}, b \neq e_k, |b| = |e_i|, 1 \leq i \leq k-1\}$, $k \leq n-1$.
$\{e_1,...,e_k\} \in \{1n,...,(n-1)n\}$, for $|e_1| < |e_2| < ... < |e_k|$.

For example, n=2 leads to the sets $T^{e1}$, $T^{e1,e2}$, $T^{e1,e2,e3}$, $T^{e1,e3}$, $T^{e2}$, $T^{e2,e3}$, $T^{e3}$. Each of these sets gives rise to distinct trees. Modified Mayeda-Seshu algorithm clearly generates all graceful trees of order n+1. However, this does not ensure the exhausting all trees. It is required to show at each level for each discarded edge as incoming edge leads to a tree isomorphic to one corresponding to some tree arising out of an incoming edge. No such easy way so far has been is known.

**Algorithm 6.2.** This has its basis partly in the following procedure: Set of node labels is $\{0,1,2,...,n\}$. Consider the following array of pairs of the type (i,i+1), i=1,2,...,n-1. The pair (i,j) is written as ij for convenience.

$S_0$ : 00   01   02   ...   0(n-1)   0n
$S_1$ : 10   11   12   ...   1(n-1)   1n
       .
       .
       .
$S_{n-1}$: (n-1)0   (n-1)1   (n-1)2   ...   (n-1)(n-1)   (n-1)n

where ij stands for the edge with end labels i and j. Define, $m = \max\{|ij|\}$, $0 \leq i,j \leq n$, $j \neq i$), or $m = \max\{|i0|, |in|\}$, $i = 0,1,...,n-1$, where $|ij|$ stands for $|i-j|$. Consider the vector $(m_1,m_2,...,m_n)$. Denote the set of all permutations $\pi(n)$ so that $a=(a_1,a_2,...,a_n) \in \pi(n)$ if and only if $a_i \leq m_i$, for i=1,2,...,n and $a_i$'s are distinct. For each $a=(a_1,a_2,...,a_n) \in \pi(n)$, get all (n+1)-tuples $g=(g_0,g_2,...,g_n)$, where $g_0= 0n$, $g_i \in A_i$ and $A_i = \{ij \in S_i : |ij| = a_i\}$, i=1,2,...,n. Designate by A(a) the set of all such (n+1)-tuples. Each such g is a graceful labeling But, the graph of it may not be a tree. This procedure exhausts all graceful labeling of trees (Prove?). How to eliminate (characterize) g's leading to a graph (disconnected with unicyclic component) which is not a tree? This happens for $n \geq 6$. (Is this why, the problem is now classified as a Graph Theoretical disease?).

All (p,q)-graphs with all graceful labelings are embedded in $K_{q+1}$. In particular, all trees of order p are embedded in $K_p$. Since edge label i occurs p-i times in $K_p$ and every choice of p-1 edges labeled 1,2,...,p-1 indices a tree and is in fact a graceful tree. Therefore,

**Proposition 6.3.** There are (p-1)! graceful labelings of all trees of order p.

**Observation 6.4.** Nodes (endnodes) of a tree may be i-repelling for some i=0,...,q. That is all endnodes of a tree need not permit label i at endnodes as well as at interior nodes. Smallest such tree is of order 6. 0-repelling nodes are shown as thick nodes for trees in **Fig. 7a**. Note that tree of order 6 is induced in the other. Subsequently, if node u is 0- or p-repelling then a caterpillar C mounted on it results in a tree with each extreme endnode away from u is also 0- or p-repelling. Therefore there exists a tree of order p for each p≥6 with the property. However, the same conclusion may not be drawn for the intermediate nodes of C. In other words, from a graph with a 0- or p-repelling node, an infinite family of graphs with this property can be constructed for every order larger than that of the given graph. Nodes repelling other node labels like 1 also do exist.

Edges, which are 0- or p-repelling at its end nodes, do exist. Some examples are shown in **Fig.**

**7b** with thick edge. Smallest tree with the property is of order 4 and is $P_3$. As in i-repelling nodes extension of a known tree with a i-repelling edge is not possible. Further, note that a tree of order 5 obtained from $P_4$ with one endnode added at the central node becomes 0- or 4-attractive.

**Problem 6.5.** Describe all graceful labelings of the class in **Fig. 7c** especially ones with 0 or n as shown**.**

In particular, consider a graceful labeling with node labels as shown in **Fig. 7c.** A label assigned to a node of B is itself an edge label. But $(n+3-\alpha)$ must be assigned to A yielding edge label $(n+3-2*\alpha)$. But then $n+3-2*\alpha = \alpha \Rightarrow \alpha = (n+3)/3$. (*Huang, Kotzig and Rosa (1982)* have studied a generalized class and classified them with respect to $\alpha$- and $\beta$- valuations.) An example of such tree is shown in **Fig. 7c**. Some other simple cases with $|A| = 2$ or 3 or $|B|$, $|A| = 1$ or $|B|$, $|A|$ prime can be handled. It will be interesting to investigate the following:

**Conjecture 6.6.** There exists a 0-attractive endnode in any tree.

In fact this can be further extended in the following way: There exists a 0-attractive noncentral endnode in any tree. Truth implies that there exists q-attractive endedge in any tree.

**Approaches.**

1.  General rule to obtain a graceful labeling of a given tree.
2.  Enlarging the known graceful trees to attain the total class of trees.
3.  Modified Mayeda-Seshu algorithm generates all possible graceful labelings of trees of order p filtering out the other labelings. But the fact that *all trees of order p are indeed encountered* in this process remains to be verified. Alternatively, a mechanism to assert that each labeled tree rejected during the generation process is indeed isomorphic to one among the gracefully labeled trees is lacking.
4.  A methodology that modifies a known graceful labeling of a tree to a graceful labeling satisfying certain properties is desired. Such results aid in obtaining a proof by Induction. Now it is apparent from literature the complexity of such a proof. This is mainly due to lack of knowledge to find a graceful labeling from a known labeling as some times no graceful labeling exists assigning a specific node label to a specified node.
5.  An algorithm which does generate a graceful labeling satisfying certain conditions from a given graceful labeling or extend a known graceful labeling of a tree to a graceful labeling of the given tree, etc.

# 7. SUPERGRACEFUL GRAPHS

The following results are proved in *Hebbare (1980).*

**Theorem 7.1.** $K_1+G$ is graceful if G is supergraceful.

**Theorem 7.2.** $K_p$ is non-supergraceful for $p \geq 3$.

This follows as a consequence of Theorem 7.1. An alternative proof of Theorem 7.2, that is enumerative follows:

**Proof.** For p≥3 the graph $K_p$ has $^pC_2$ edges. Let $r = {}^{p+1}C_2$. If $K_p$ were supergraceful, assign a subset N of p labels from S={1,2,...,r} to the nodes in such a way that the edges receive each of the labels from {1,2,...,r}-N. Let φ be a total labeling of $K_p$. Note that for any label $s \in S$ either $s \in N(\varphi)$ or $s \in E(\varphi)$. The label $r \in N(\varphi)$ only. Next, if $r-1 \in N(\varphi)$ then $1 \in N(\varphi)$. Assert $r-2 \notin N(\varphi)$. Otherwise 1 repeats as edge label. Hence $2 \in N$ but then r-2, $r-3 \in E(\varphi)$. Lastly, $r-4 \notin N(\varphi)$ ((r-4) $\notin E(\varphi)$) otherwise r-3 repeats. If $r-1 \in E(\varphi)$ then $q \in N(\varphi)$. Now, $r-2 \notin E(\varphi)$ otherwise r-1 repeats. If $r-2 \in E(\varphi)$ then $2 \in N(\varphi)$ in which case 1 repeats as a node and edge labels. Thus if p>3 a total labeling is impossible. ☐

Another interesting consequence of Theorem 7.1 is given in the following:

**Theorem 7.3.** If $H = K_1 + G$ and H is non-graceful then G is non-supergraceful.

**Corollary 7.4.** If $H = K_1 + G$ is an eulerian graph with $q \equiv 1 \text{ or } 2 \pmod 4$ then G is non-supergraceful.

**Example 7.5.** $K_n - e$, where e is an edge of $K_n$, is known to be non-graceful for n=6 and 10. Hence by Theorem 1, $K_5 - e$, $K_9 - e$ are non-supergraceful.

Some classes of non-supergraceful graphs arising from Corollary 2 are:

**Example 7.6.** The wheels $W_n$, $n \equiv 3 \pmod 4$ are non-supergraceful since $K_1 + W_n$ is an eulerian graph with q edge $q \equiv 2 \pmod 4$. Similarly the wheels $W_n$, $n \equiv 2 \pmod 4$ with a new node attached at the center of the wheel is non-supergraceful.

**Example 7.7.** The n-prism $P_n = K_2 \times C_n$, n>3 is a (2n,3n)-cube graph. It is easy to check that $K_1 + \bar{N}?(n)$ is an eulerian graph with q edges, $q \equiv 1 \text{ or } 2 \pmod 4$ whenever $n \equiv 1 \text{ or } 2 \pmod 4$ and hence non-supergraceful.

**Example 7.8.** In fact, any regular graph G of odd degree such that $K_1 + G$ is an eulerian graph with q edges, $q \equiv 1 \text{ or } 2 \pmod 4$ is a non-supergraceful graph. Petersen graph and Dodecahedron are examples of such non-supergraceful graphs.

Should a regular graph G of order p and degree d with $p \equiv 1 \text{ or } 3 \pmod 4$ and $d \equiv 2 \pmod 4$ the graph $H = K_1 + G$ is a non-supergraceful graph? Among the five platonic solids the cube, $Q_3$ is supergraceful as shown in **Fig. 8a** and all the others are non-supergraceful. Tetrahedron are ruled out by Corollary 7.4. The Octahedron is non-supergraceful since Octahedron with a new node adjacent to each of its nodes is non-graceful.

**Theorem 7.9.** RRK- tree conjecture implies that trees are supergraceful

**Proof.** If a tree T is graceful and φ is a graceful labeling of T, then $N(\varphi) = \{0,..., p-1\}$. We now exhibit a total labeling of T. Define a new labeling μ of T as follows: $\mu(x) = \varphi(x) + p$, for each $x \varepsilon V$. Then $N(\mu) = \{p, p+1,..., p+q\}$ and $E(\mu) = E(\varphi)$. Hence μ is a total labeling and T is supergraceful. ☐

We remark that not every total labeling of T is of type in the proof of the above theorem. But if such a labeling exists for a tree then it is graceful. That '*trees are supergraceful*' was conjectured by *Acharya [2]*. Graphs with p≤5 nodes are supergraceful except for the three graphs $K_4$, $K_5$-e and $K_5$. This is verified by the total labelings given in **Fig. 9.** At least eight non-supergraceful graphs of order six are known. If G is a unicyclic graph such that $K_1+G$ is eulerian then it has q edges, $q \equiv 0 \pmod 4$. Cycles $C_n$, $n \geq 3$ are shown supergraceful in *Acharya (1983)*. That unicyclic graphs of order at most six are supergraceful has been verified and unicyclic graphs of order 6 are shown with a total labeling in **Fig. 10**. We now make the following:

**Conjecture 7.10.** Unicyclic graphs are supergraceful.

## 8. SUPERGRACEFUL EMBEDDINGS

Before considering the supergraceful embeddings we shall digress to establish some results for (non) supergraceful graphs.

**Observation 8.1.** Let G be a non-supergraceful graph and $\varphi$ be a semitotal labeling of G. Then the labeling $\mu$ of $H=K_1+G$ such that $\mu_V=\varphi$ and $\mu(u)=0$, where $K_1=\{u\}$ is an optimal labeling of H and conversely.

It is conjectured by *Acharya [2]* that for every (p,q)-graph G there exists a semitotal labeling such that $1 \in N(\varphi)$. This is true for a special class of graphs as follows:

**Theorem 8.2.** Let G be a graph with a full degree node (that is, a node of degree p-1) u and there is a semitotal labeling $\varphi$ such that $\varphi(u) = opt(G)$ and $opt(G)-1 \in N(\varphi) \cup E(\varphi)$, where $opt(G) = p+q+q_0$, $q_0>0$ is an integer. Then there is a semitotal labeling $\mu$ of G such that $1 \in N(\mu)$.

**Proof.** Firstly note that $opt(G)$ must be a node label for any semitotal labeling of G. If $1 \in N(\varphi)$ then nothing to prove. Otherwise by assumption $opt(G)-1 \in N(\varphi)$. Define a new labeling $\mu$ as follows: $\mu(u) = \varphi(u)$, where $\varphi(u) = opt(G)$ and $\mu(x) = opt(G)-\varphi(x)$, for each $x \in V-\{u\}$. Clearly, $\mu$ is a semitotal labeling as the edge labels in G-u are unaltered and the node (edge) labels of the nodes $x \in V-\{u\}$ (the edge labels of e= ux) are interchanged. Now since $opt(G)-1 \in N(\varphi)$ we conclude that $1 \in N(\varphi)$. ∎

**Corollary 8.3.** If G is a supergraceful graph with a full degree node u and a total labeling such that $\varphi(u) = p+q$ then there exists a total labeling $\mu$ of G such that $1 \in N(\varphi)$.

**Observation 8.4.** Let G be a $(p_1,q)$-graph. Then there exists a graceful $(p,q+(p-p_1))$-graph for each $p>p_1$ containing G as an induced graph. Introducing new nodes adjacent to the node u labeled zero, assigning node labels $q+1, q+2,..., q+(p-p1)$, can do this.

The following theorem describes a similar construction for supergraceful graphs.

**Theorem 8.5.** Let G be a supergraceful (p,q)- graphs with a total labeling such that $1 \in N(\varphi)$. Then there exists a supergraceful (p+t,q+t)-graph for each $t \geq 0$.

**Proof.** Define a new a (p+t,q+t)-graph H containing G as follows: Let u be such that $\varphi(u)=1$. Attach $t \geq 0$ nodes at u and define $\mu$ as below: $\mu_V=\varphi$ and $\{\mu(x): x \in V(H)-V\} = \{p+q+2, p+q+4,..., p+q+2t\}$. $\mu$ can be verified to be a total labeling of H. ∎

The above theorem can be further extended as follows:

**Observation 8.6.** Let G be a supergraceful (p,q)-graph and φ be a total labeling of G such that α = min{φ(x):x∈V}, α≥1. Then there exists a supergraceful graph of order p+αt for each t≥0, containing G as an induced subgraph.

Following was conjectured in *Acharya [2]*. We exhibit a disconnected embedding and consider connect embedding whenever G has a semitotal labeling with 1 as a node label:

**Theorem 8.7.** Every graph can be embedded into a supergraceful graph.

**Proof.** Let G be a (p,q)-graph. If G itself is supergraceful then nothing to prove. Otherwise G is non-supergraceful. Define a graph H containing G and a isolated nodes with weights assigned forming the set $\{1,2...,p+q+\alpha\}-(N(\varphi) \cup \overline{N}(\varphi))$ where φ is a semitotal labeling of G. □

The graph H in the above theorem need not be connected embedding of G with G as an induced subgraph. Does there exist a connected supergraceful embedding for a graph? The following theorem asserts this whenever G has a semitotal labeling φ such that $1 \in N(\varphi)$.

**Theorem 8.8.** Let G be a graph with a semitotal labeling φ such that $1 \in N(\varphi)$. Then G can be embedded as an induced subgraph into a connected supergraceful graph.

**Proof.** By definition G is non-supergraceful and let φ be a semitotal labeling of G with $1 \in N(\varphi)$ and φ(u)=1. Any $x \in \overline{N}(\varphi) \cup \overline{E}(\varphi)$ may occur in the embedding as a node or edge label. Two cases arise according as opt(G) is even or add. We shall prove the theorem in the second case and the first case can be proved on the similar lines. Define, $E_1 = \{x \in \overline{N}(\varphi) \cup \overline{E}(\varphi): x \text{ is even}\}$, $E_2 = \{x \in \overline{N}(\varphi) \cup \overline{E}(\varphi): x \text{ is odd}\}$, $m_1 = \min E_1$, $m_3 = \max E_1$, $m_2 = \min E_2$, $m_4 = \max E_2$. Select $v,w \in V$ such that $\varphi(v) = \min \{\varphi(x): \varphi(x)+m_1 > d(G), \varphi(x) \text{ is even}\}$, $\varphi(w) = \min \{\varphi(x): \varphi(x)+m_2 > opt(G), \varphi(x) \text{ is odd}\}$. Further, $t = \max \{\varphi(v)+m_3, \varphi(w)+m_4\}$.

Now, construct a graph H containing G with t additional nodes adjacent to u and labeled by opt(G)+2, opt(G)+4,...,t. Now, even (odd) $x \in \overline{N}(\varphi) \cup \overline{E}(\varphi)$ can be obtained by joining appropriate new node to φ(v) (φ(w)). The resulting labeling of H can be verified to be a total labeling and hence H is a supergraceful graph. □

## 9. PRIME GRAPHS

A `*node labeling*' α of a graph G of order n is a bijection α: V(G) ⟶ {1,2,...,n }= $N_n$. α is said to be a `*prime labeling*' of G if for each edge uv ε E(G), (α(u), α(v))=1, that is, α(u) and α(v) are coprime to each other. A graph that admits a prime labeling is called a `*prime graph*'. Henceforth, the node labels of a graph G are taken to be the values assigned to the nodes by a prime labeling α, when we are dealing with G and α. The work on prime graphs is initiated in *Tout, Dabboucy and Howalla [5],* although the concept is due to R.C. Entringer. Further, he made the following:

**Conjecture 9.1.** *(Entringer)*. All trees are prime.

Truth of the conjecture implies that all the forests are prime. This follows from the hereditary property of prime graphs (see Property 4 below). Some families of trees including caterpillars with constant number of pendant nodes adjacent to its internal nodes and binary trees were shown to be prime in [1,5] which support the Conjecture. The result on caterpillars is now extended to any prime graph in *Acharya (1983)*. For some other infinite families of prime graphs see *[1,5]*.

Define `*super prime graph'* $SP_n$ of order n as follows: $V(SP_n)=\{1,2,...,n\}$ and two nodes uv ε $V(SP_n)$ are adjacent if and only if $(u,v)=1$. $SP_n$ is maximal and unique and hence the name. Here, we list some interesting properties of $SP_n$ graph and prove some graph theoretical properties. Note that the study of structural properties of $SP_n$s is significant in the study of prime graphs.

**Elementary Properties.** Number of edges m in $SP_n$ is given by $m=\sum \varphi(x)$, x=1,..,n-1 where $\varphi(1)=0$, $\varphi(x)=\{t<x: (x,t)=1\}$. $\varphi$ is the Euler's $\varphi$-function and so is a necessary condition for an (n,m)-graph to be prime, that is an (n,m')-graph is not prime for m' > m. Table below gives values of m for n≤20.

| Table 2. | n: | 1 | 2 | 3 | 4 | 5 | 6 | 7 | 8 | 9 | 10 | 11 | 12 | 13 | 14 | 15 | 16 | 17 | 18 | 19 | 20 |
|---|---|---|---|---|---|---|---|---|---|---|---|---|---|---|---|---|---|---|---|---|---|
| | m: | 0 | 1 | 3 | 5 | 9 | 11 | 17 | 21 | 27 | 31 | 41 | 45 | 57 | 63 | 81 | 89 | 105 | 111 | 129 | 137 |

Let $\mathcal{P}_n = \{x \; \varepsilon \; N_n : x < n \text{ and } x \text{ is prime}\}$. Then the size of the largest clique in $SP_n$ is given by $\omega(SP_n) = |\mathcal{P}_n \cup 1|$, and in fact $<\mathcal{P}_n \cup 1>$- the graph induced by $\mathcal{P}_n \cup 1$ in $SP_n$, is the largest clique. Further, the independence label is given by $\beta_o(SP_n) = [n/2]$. In fact, the set of all even labeled nodes E is the maximum independent set. Hence, for any prime graph G the following are necessary conditions. $\omega(G) \leq \mathcal{P}_n+1$, and $\beta_o(G)>|n|$ (See also *[1]*). Lastly, note that any clique of size at least two in a prime graph G yields a set of mutually coprime labels.

Let $\underline{\mathcal{P}}_n = \{x \varepsilon N_n: x \; \varepsilon \mathcal{P}_n \text{ and } 2x < n\}$. Then $SP_n$ has exactly $|\underline{\mathcal{P}}_n|+1$ nodes of degree n-1 and this set is precisely $\underline{\mathcal{P}}_n \cup 1$.

Every spanning subgraph of a prime graph is prime. In particular, every spanning subgraph of $SP_n$ is a prime graph and conversely. Further, every prime graph of order ≤n is a subgraph of $SP_n$. This justifies name the super prime graph. Hereditary property of prime graphs of $SP_n$ is significant. The problem of finding all prime graphs of order n reduces to finding all the distinct spanning subgraphs of $SP_n$. Lastly, if G is not a prime graph of order n then every supergraph of G of order n) is not a prime graph.

A graph G of order n is said to be `*tree-complete'* if every tree of order n is a subgraph of G (see *Nebesky [5]*). In view of the above definition and the property 4, Entringer's Conjecture is equivalent to the following:

**Conjecture 9.2.** Super prime graph is tree-complete.

Let G be a prime graph of order n. Then G-n is prime. In particular, if n is a pendant node of a prime tree T then T-n is again a prime tree of order n-1. We now propose the following strong prime tree Conjecture:

**Conjecture 9.3.** For every tree of order n there exists a prime labeling which assigns n to a pendant node.

**Observation 9.4.** If G is a tree-complete graph then $K_1+G$ is tree-complete where `+' denotes the sum of graphs. As a consequence of Theorem 3 in *[4]* and the Property 3 we get the following:

**Theorem 9.5.** Let $i\varepsilon\{1,2\}$ and let $SP_n$ be such that every tree $T_o$ of order n with $\Delta(T_o)<[(p+i)/2]$ is isomorphic to a spanning subgraph of G. Then $SP_{n+1}$ and $SP_{n+2}$ are tree complete. Further, if n+1 and n+3 are twin primes then $SP_{n+1},\ldots,SP_{n+5}$ are tree-complete.

As a consequence of the above theorem and Observation 1 we have the following:

**Corollary 9.6.** Let $SP_n$ be a tree-complete graph such that n+1 is prime. Then $SP_{n+1}$ and $SP_{n+2}$ are tree complete. Further, if n+1 and n+3 are twin-primes then $SP_{n+1},\ldots,SP_{n+5}$ are tree-complete.
Using Corollary 1 in [4] and the Corollary 2, Conjecture 2 is verified for all $n\leq 15$. Note that, if $SP_{16}$ is proved to be tree-complete then since 17 and 19 are twin-primes, it follows from Corollary 2 that SP(x) is tree-complete for x=17, 18, 19, 20, 21.

Define $\mathcal{P}_{m,n}$ to be the product of primes p, $m\leq p\leq n$ for $m,n \varepsilon \mathcal{P}_n$. When m=n=p then $\mathcal{P}_{m,n}=\mathcal{P}_p$. Then the node of $SP_n$ labeled $\mathcal{P}_{2,t}$ such $\mathcal{P}_{2,t}<n$, for t=max $\mathcal{P}_n$, has minimum degree and in fact, $\delta(SP_n)=d(\mathcal{P}_{2,m})= n-[n/2])$. Yet this is another necessary condition implying that any graph G of order n with $\delta(G)>\delta(SP_n)$ is not prime. $\mathcal{P}_{2,p}$ is even for any $p\varepsilon\mathcal{P}_n$. $d(p)=n-[n/p]$. Let r be the number of distinct degrees in $SP_n$. Further for $p\varepsilon\mathcal{P}_n$, $N(p)=N(p^2)=\ldots=N(p^a)$, where N(x) denotes the neighborhood of x in $SP_n$ and a is maximum such that $p^a<n$. Also, for $p_1,p_2\varepsilon\mathcal{P}_n$, $N(p_1p_2) = N(p_1)\cap N(p_2)$. Combining the above two results, we get for $n_1=p_1^{\alpha_1}p_2^{\alpha_2}\ldots p_c^{\alpha_c}$, that $N(n_1)=N(p_1^{\alpha_1})\cap\ldots\cap N(p_c^{\alpha_c})$, where $2\leq n_1\leq n$. Designate the set of distinct primes in the prime factorization of n by $d\mathcal{P}(n)$. For $2\leq n_1,n_2\leq n$, such that one of $n_1, n_2$ is not prime $d\mathcal{P}(n_1)=d\mathcal{P}(n_2)$ implies that $d(n_1)=d(n_2)$. Each node of $\mathcal{P}_n\cup 1$ has full node, that is of degree n-1. $\mathcal{P}_{>2}=\{x:2x>n\}$ and $\mathcal{P}_{<2}=\{x:n/2\leq x\leq n\}$. For each prime $p\in\mathcal{P}-\mathcal{P}_n$, nodes labeled p, p2,…,$p^\alpha$ have same degree where $p^\alpha<n$ and $\alpha$ is maximum with this property. Lastly, let $\tau$ be the possible distinct products of the type $p_1p_2\ldots p_\alpha$ of $\alpha>1$ primes with the product <n. Now we can enumerate distinct degrees in $SP_n$:

**Observation.** Distinct degrees in $SP_n$ is given by $|\mathcal{P}_n\cup 1|+|\mathcal{P}_{>2}|-|\mathcal{P}_{<2}|+\tau$.

**Observation 9.7.** $K_{r,s}$ is prime if and only if r or $s<|\mathcal{P}_n|+1$, where $\mathcal{P}_n$ is defined as in Property 3 (See also [1]). Further, when r or s = $|\mathcal{P}_n|+1$ the graph has maximum number of edges.

## 10. EXTREMAL NON-PRIME GRAPHS

**Observation 10.1.** *([1,4])* Every prime graph G of order n satisfies the following inequality. $\beta_o(G)\geq[n/2]$ where $\beta_o(G)$ is the independence number of G.

Paul Erdos asked the following question at the third MATSCIENCE Conference on Label Theory (3-6 June,1981) at Mysore.

**Problem 10.2.** What is the minimum number of lines a non-prime graph of order n should have?

We consider this in two cases according as the graph is connected or not. Let $\mu(n)$ and $\mu'(n)$

denote the minimum number of lines a non-prime graph can have in the respective cases. Below we obtain upper bounds for $\mu(n)$ and $\mu'(n)$ and conjecture that the bounds should be attained.

**Connected Case.** Note that $K_4$ and $K_5$ are the only non-prime graphs of order at most five. Hence, $\mu(n)$, $n \leq 3$ is not defined, $\mu(\text{mod } 4) = 6$ and $\mu(5) = 10$. In the case n=6, $\mu(6) = 7$ and there is a unique (6,7)-graph which is not prime as shown in **Fig. 12a.** The graph in **Fig. 12b** shows that $\mu(7) \leq 10$ and conjecture that $\mu(7) = 10$ holds. We prove the following:

**Theorem 10.3.** If $n \geq 8$ then $\mu(n) \leq \begin{cases} n+1, \text{ if n is even,} \\ n+2, \text{ otherwise.} \end{cases}$

**Proof.** Inequalities are established by exhibiting a suitable non-prime graph.

**Case 1.** n even. Let n=2t+6, t≥1. Consider the (n,n+1)-graph $G_1$ in **Fig. 12c**. Clearly, $\beta_o(G_1) \geq t+2$ since $\{v_1,...,v_t, u, v\}$ is an independent set, where u and v are one node from each of the two triangles. On the other hand if A is an independent set of $G_1$ then A can have at most one node from each of the two triangles and $|A \cap \{u_1,...,u_t, v_1,...,v_t\}| < t$. Hence, $|A| < t+2$. Thus, $\beta_o(G_1) = t+2$ follows. But then by Observation ?, $G_1$ is not prime.

**Case 2.** When n is odd, we consider two cases according as n≡1 or 3(mod 4). In both these cases the graphs $G_2$, $G_3$ given in **Fig. 12d,e** can be verified to be non-prime graphs. *Case 2(a).* n ≡ 1(mod 4). Consider the following (4t+1,4t+3)-graph, for t>1. One can check that, $\beta_o(G_2) = 2t-1$. *Case 2(b).* n≡3(mod 4). Consider the following (4t+3, 4t+5)-graph for each t>1. One can check that $\beta_o(G_3)=2t$.

**Conjecture 10.4.** For a connected graph of order n (≥8), equality in Theorem10 holds for $\mu(n)$.

Truth of the conjecture implies truth of Entringer's prime tree conjecture and further implies that *all unicyclic graphs are prime*.

**Disconnected Case.** Let G be a graph not necessarily connected. One can verify that all disconnected graphs with at most five nodes are prime and hence $\mu'(n)$, n≤3 is undefined and $\mu'(n) = \mu(n)$ for n=4,5. Further, $\mu'(6)=6$ and there is a unique disconnected non-prime graph namely, $2K_3$, that is, two copies of $K_3$. That $\mu'(7) \leq 9$ follows from the non-prime graph $K_3 \cup K_4$, and conjecture that $\mu'(7)=9$ holds. We prove the following:

**Theorem 10.5.** If n≥8 then $\mu'(n) \leq n$.

Proof. Graphs, $tK_3$, if n=3t, t≥3; $(t-1)K_3 \cup C_4$, if n=3t+1, t≥3 and $(t-1)K_3 \cup C_5$, if n=3t+2, t≥2 can be verified to be non-prime using observation and satisfy the bound in the theorem.

**Conjecture 10.6.** For any graph of order n, n≥8, $\mu'(n)=n$ holds.

Lastly, we propose the following for regular graphs:

**Conjecture 10.7.** 0% of all (connected) regular graphs are not prime.

This conjecture is supported by the result in [1] that every non-bipartite regular graph G with $\beta_o(G) = [n/2]$ is not prime. Further, it is known that $K_{r,r}$ is prime whenever n≤2 and $K_n$ is prime whenever n≤3. A regular graph whose components are odd cycles is clearly not prime. In fact, a regular graph of degree 2 with sufficiently many odd cycle components can be shown to be a non-prime graph. On the other hand, cycles are prime and the n-dimensional cube $Q_n$ is known to be prime for n≤3 only. We remark that $Q_3$ is the only prime graph out of the five Platonic solids. Although $K_{3,3}$ is not prime, $2K_{3,3}$ is a prime graph. A scattered example of a regular prime bipartite graph is given in **Fig. 12f.**

The five e-minimal non-prime graphs of order six are $2K_3$, $K_{3,3}$, $K_5 \cup K_1$ and the two graphs of **Fig. 12g,h**. Note that, a nonprime graph G with µ'(G) edges is minimal but not conversely.

## Acknowledgements

Author is gratefull to S/Shri. **T.K.N. Gopalaswamy,** Former Director (Exploration), **S.K. Patra,** G.G.M. (Exploration), **A.V. Raju,** G.M. (Reservoir) and **S.K. Majumder,** Head M.S.G.(E) and G.M.(Geology) of Oil and Natural Gas Corporation Ltd., for permitting to attend and deliver lectures for '*Group Discussion on Graph Labeling Program*' at KMREC, Surathkal (16-25 August 1999) and publish this work with necessary facilities. Thanks are also due to Dr. **M.S. Rao**, Dy. G.M. (Programming) and **G.C. Raturi,** Chief Manager (P&A) of Oil and Natural Gas Corporation Ltd., for computing and administrative facilities.

Author is thankful to Dr. **B.D. Acharya,** (author's friend and colleague since 1976) Scientist G, Department of Science and Technology and Dr. **S.M. Hegde,** Professor, Mathematics Department, KREC, Surathkal for extending the invitation to participate in *'Group Discussion on Graph Labeling Problems'* and for hospitality.

## REFERENCES

1. **B.D. Acharya,** *On d-Sequential Graphs,* J. Mathematical Physics Sciences, 17(1983) 21-35.
2. **B.D. Acharya,** *Embedding Graphs in Graceful graphs*, Report (1980), MRI, Allahabad.
3. **B.D. Acharya,** *Supergraceful Graphs.* Report (1980), MRI, Allahabad.
4. **B.D. Acharya,** *Construction of Certain Infinite Families of Graceful Graphs from a Given Graceful Graph*. Report (1980), MRI, Allahabad.
5. **J. Ayel and O. Favaron**, *The Helms are Graceful*, Progress in Graph Theory, Academic Press, Toronto 89-92 (1984).
6. **D.W. Bange and A.E. Barkauskas and P.J. Slater,** *Simply Sequential Graphs*, Proceedings of 10th S-E Conference Combinatorics, Graph Theory and Computing, (1979) 155-162.
7. **J.C. Bermond,** *Graceful Graphs, Radio Antenne and French Windmills,* Proceedings of Open University Conference, May 1978, 18-37.
8. **V. Bhat-Nayak and S.K. Gokhale,** *Validity of Hebbare's Conjecture*, Utilitas Mathematica, 29(1986) 49-59.
9. **G.S. Bloom,** *Dissertation*, University of Southern California, 16-17, (1975).
10. **G.S. Bloom**, *A Chronology of the Ringe-Kotzig conjecture and the Continuing Quest to call all the Trees Graceful.* **'Topics in Graph Theory'** Ed. F. Harary. New York Academy of Sciences (1979) 32-51.
11. **G.S. Bloom and S.W. Golomb,** Applications *of Numbered Undirected* Graphs, Proceedings IEEE 65 (1977) 562-570.
12. **G.S. Bloom and S.W. Golomb,** *Numbered Complete Graphs, Unusual Rulers, and Assorted Applications,* Springer Lecture Notes on Mathematics, 642(1978) 53-65.
13. **F. Buckly and F. Harary (1988),** *Distances in Graphs'*, Addison-Wesley.
14. **G.J. Chang, D.F. Hsu and D.G. Rogers,** *Additive Variation of a Graceful Theme: Some Results on Harmonius and Other Related Graphs,* Congress Numerentium, 32 (1981)181-197.
15. **W.K. Chen, (1971),** *Applied Graph Theory*, North Holland, Amsterdam**.**

16. **Christofides N, (1975),** *Graph Theory: An Algorithmic Approach,* Academic Press.
17. **Prabir Das and T. Gangopadhyay,** *Some Results of α-Valuations and Average Numberings,* Utilitas Mathematica, 33 (1983) 209-222.
18. **C. Delorme, Maheo, H. Thuillier, K.M. Koh and H.K. Teo,** *Cycles with a Chord are Graceful*, J. Graph Theory, 4 (1980) 409-410.
    Sufficient conditions for $GxK(2)^n$ to be graceful. Average numbering is introduced.
19. **J. Doma,** *Unicyclic Graceful Graphs*, M.S. Thesis, McMaster University, 1991.
20. **E. Flandrin, I. Fournier and A. Germa,** Numerations Gracieuses Des Chemins, Ars Combinatoria, 16(1983) 149-181.
    For every path P of length ≥8, for any vertex u and any k in $Z_n$ there exists a graceful numbering which assigns k to u.
21. **D. F. Hsu,** *Harmonious Labelings of Windmill Graphs and Related Graphs*, J. Graph Theory, 6(1982) 85-87.
22. **J.G. Gallion**, *A Dynamic Survey of Graph Labeling*, Electroninc Journal of Combinatorics, 5 (1998), #DS6.
23. **T. Gangopadhyay and S.P. Rao Hebbare*,** *Bigraceful Graphs-I.* Utilitas Mathematica, 17(1980) 271-275.
24. **S.W. Golomb,** *How to Number a Graph,* **'Graph Theory and Computing'** (Ed. R.C. Read), Academic Press, New York (1972) 23-37.
25. **S.W. Golomb,** *The Largest Graceful Subgraph of the Complete Graph*, American Mathematical Monthly, 81 (1974) 499-501.
26. **R.K. Guy**, *Monthly Research Problems, 1969-77*, American Mathematical Monthly, 84 (1977) 807-815.
27. **F. Harary,** *'Graph Theory',* Addison Wesley, NewYork (1972).
28. **Harary and D.F. Hsu,** *Node-Graceful Graphs*, Mathematical Applications, 15(1988) 291-298.
29. *'Topics in Graph Theory'*, Ed **F. Harary**. NYAS, (1979) 32-57.
30. **Hebbare S.P. Rao*,** *Graceful Cycles*, Utilitas Mathematica, 10 (1976) 307-317.
31. **Hebbare S.P. Rao*,** *Graceful Algorithms and Conjectures*, Report (1980), MRI, Allahabad.
32. **Hebbare S.P. Rao*,** *Supergraceful Graphs,* Report (1980), MRI, Allahabad.
33. **Hebbare S.P. Rao*,** *Non-supergraceful Graphs*, Report (1980), MRI, Allahabad.
34. **Hebbare S.P. Rao*,** *Embeddings Graphs into Graceful and Supergraceful Graphs*, Report (1980), MRI, Allahabad.
35. **C. Huang, A. Kotzig and A. Rosa,** *Further Results on Tree Labelings*, Utilitas Mathematica, 21C(1982) 31-48.
    Labeling trees with at most four end vertices. Trees with β- labelings but without α- labelings. Labeling trees with diameter at most 4. Four vertices with at most four vertices of one color.
36. **K. M. Koh and T. Tan,** *Two Theorems on Graceful Trees.* Discrete Mathematics, 25(1979) 141-148.
37. **K.M. Koh, D.G. Rogers and C.K. Lim,** *On Graceful Graphs I: Sum of Graphs,* Sea Bull. Math. 3(1979)58.
38. **K.M. Koh, D.G. Rogers and T.Tan**, *Products of Graceful Trees*, Discrete Mathematics 31(1980) 279-292.
39. **K.M. Koh, D.G. Rogers, P.Y. Lee and C.W. Toh,** *On Graceful Graphs V: Unions of Graphs with one Vertex in Common*, Nanta Mathematica, 12(1979) 133-136.
40. **A. Kotzig,** *On Certain Vertex-Valuations of Finite Graphs,* Utilitas Mathematica, 4(1973) 261-290.
    Nearly all trees have ∝- valuations.
41. **Sin-Min Lee and Ping-Chyan Wang,** *On k-Gracefulness of the Sequential Join of Null graphs*, Pre-print.
42. **M. Maheo,** *Strongly Graceful Graphs,* Discrete Mathematics 29(1980)39-46.
    Strongly graceful, a particular case of α- valuation. If g is strongly graceful then so is G+K(2).
43. **M. Maheo and H. Thuillier,** *On d-Graceful Graphs, Ars Combinatorics*, 13(1982)181-192.
    *Generalization of graceful graphs to d- graceful graphs. Wind mills, cycles, m K(3) with K(2) in common, wheels are studied.*
44. **W. Mayeda,** (1972) *'Graph Theory'*, John Wiley.
45. **L. Nebesky,** *On Tree Complete Graphs,* Casopis pro pestovani matematiky, 100(1975)334-338.
46. **A. Rosa,** *On Certain Valuations of the Vertices of a Graph*. *'Theory of Graphs',* Garden and Breach, NY (1967) 349-355.
47. **A. Rosa**, *Labeling Snakes*, Ars Combinatoria, 3 (1977), 67-74.
48. **D.A. Sheppard,** *The Factorial Representation of Balanced Labeled Graphs,* Discrete Mathematics, 15 (1976), 379-388.
49. **P.J. Slater**, *On k-Sequential and Other Numbered Graphs*, Disc. Math., 34(1981)185-193.
50. **P.J. Slater**, *On k-Graceful, Locally Finite Graphs*, Combinatorial Theory, B 15,3,(1981)319-322.
51. **M. Turuszcznski,** *Graceful Unicyclic Graphs,* Demonstratio Mathematica, 17 (1984) 377-387.

\* Name changed from Suryaprakash Rao Hebbare to Suryaprakash Nagoji Rao, ONGC office order No. MRBC/Expln/1(339)/86-II, Dated 08 March 1999 as per Gazette notification of Govt. Maharashtra December 1998.

### Fig. 1a  Possible Higly Gracefulness

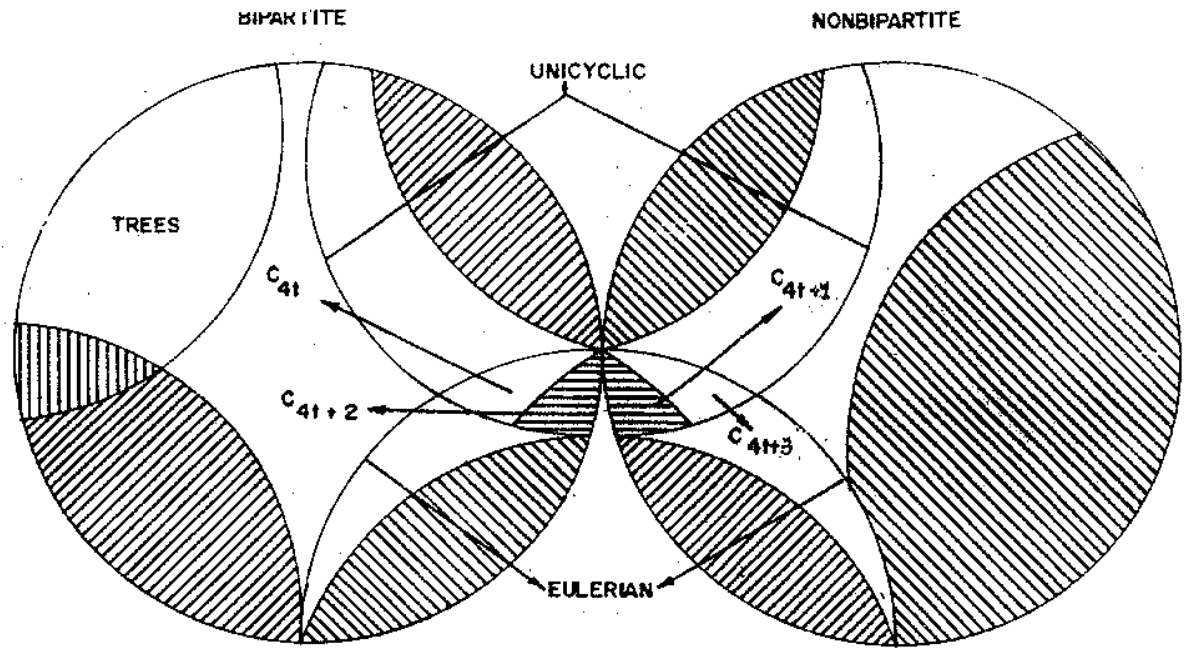

### Fig. 1b Expected Highly Gracefulness

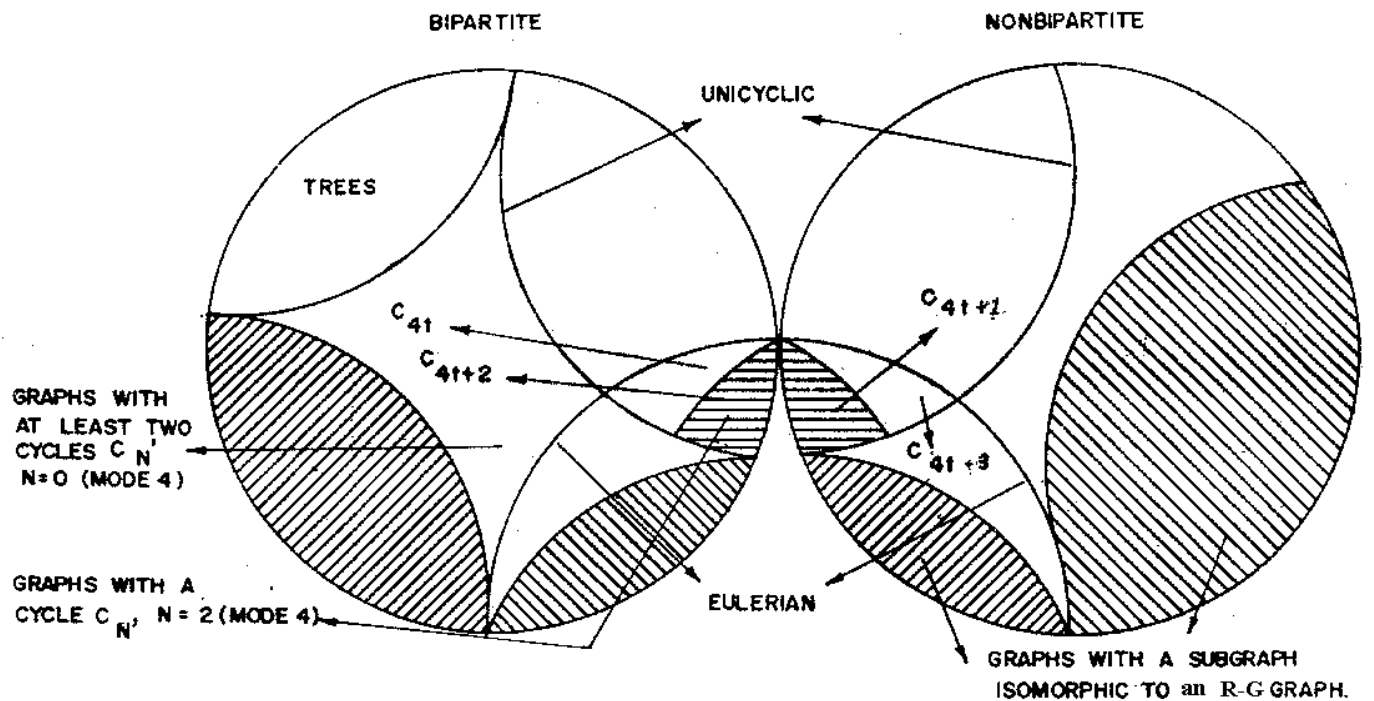

### Fig. 2a Possible Criticality

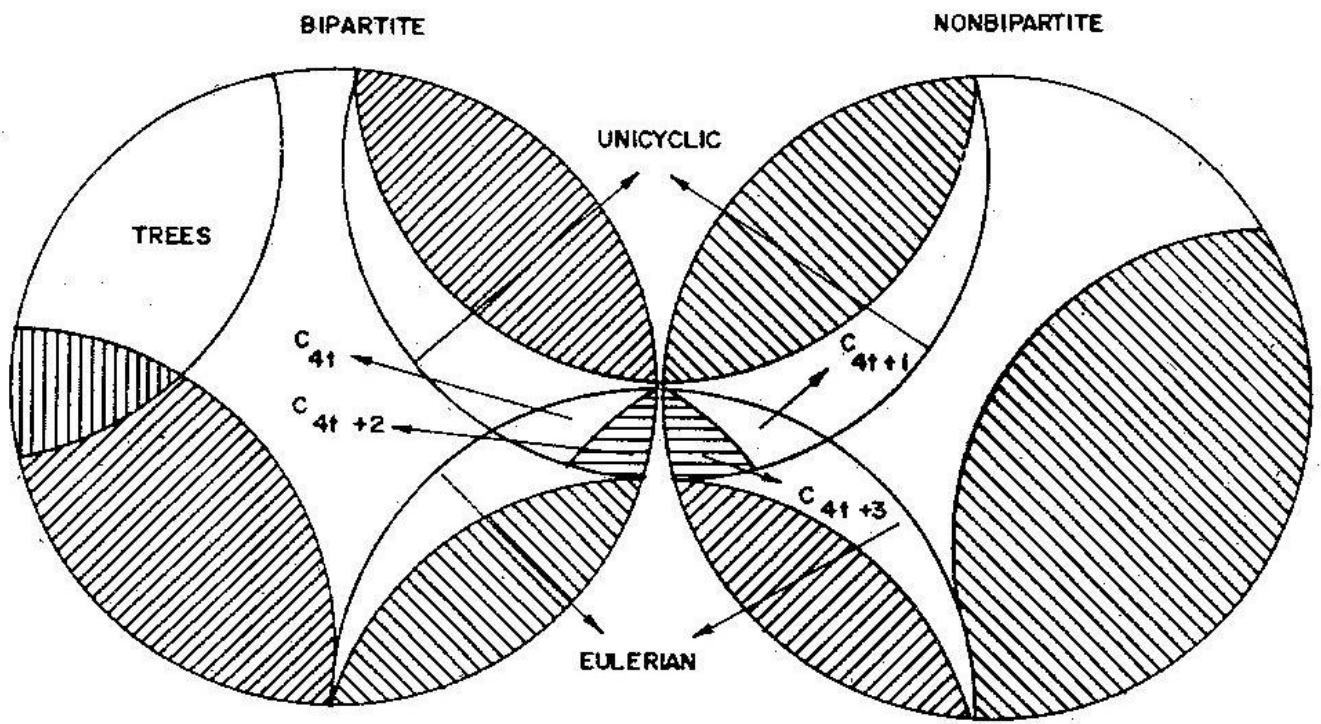

### Fig. 2b Expected Criticality

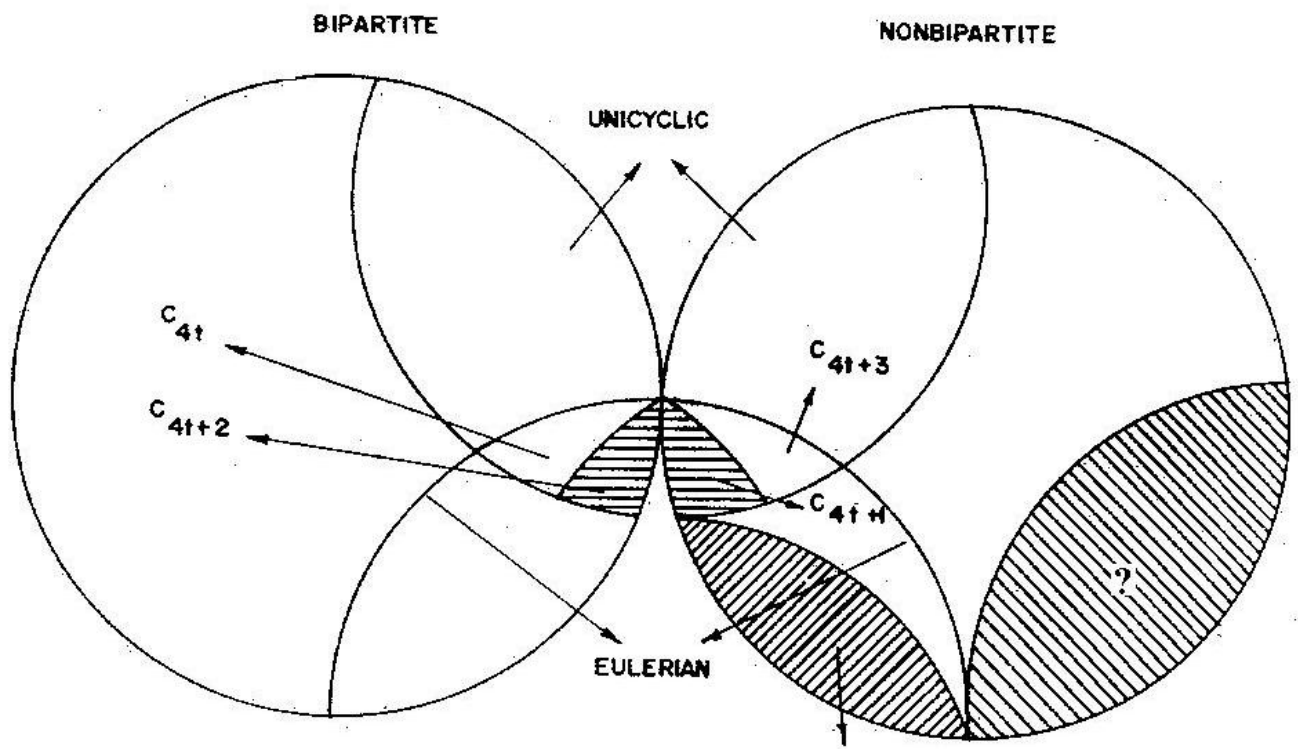

Graphs of Type 1, 2 of Proposition 3.1
and Some Others (?)

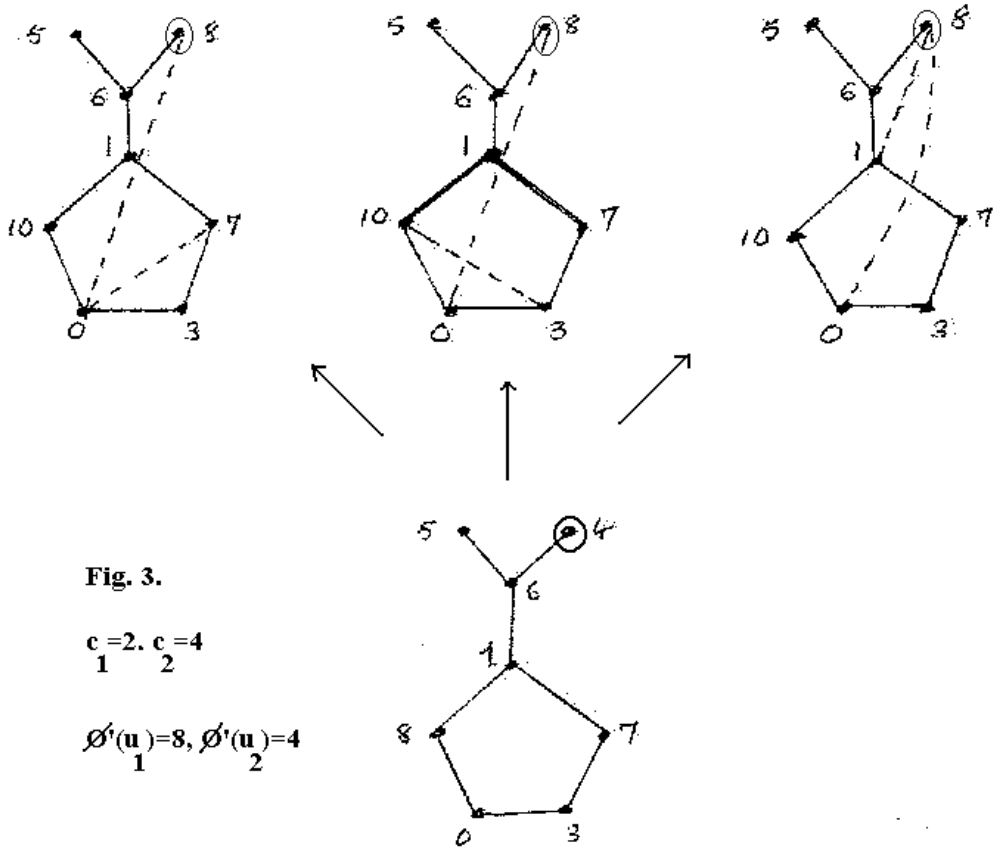

**Fig. 3.**

$c_1 = 2, c_2 = 4$

$\emptyset'(u_1) = 8, \emptyset'(u_2) = 4$

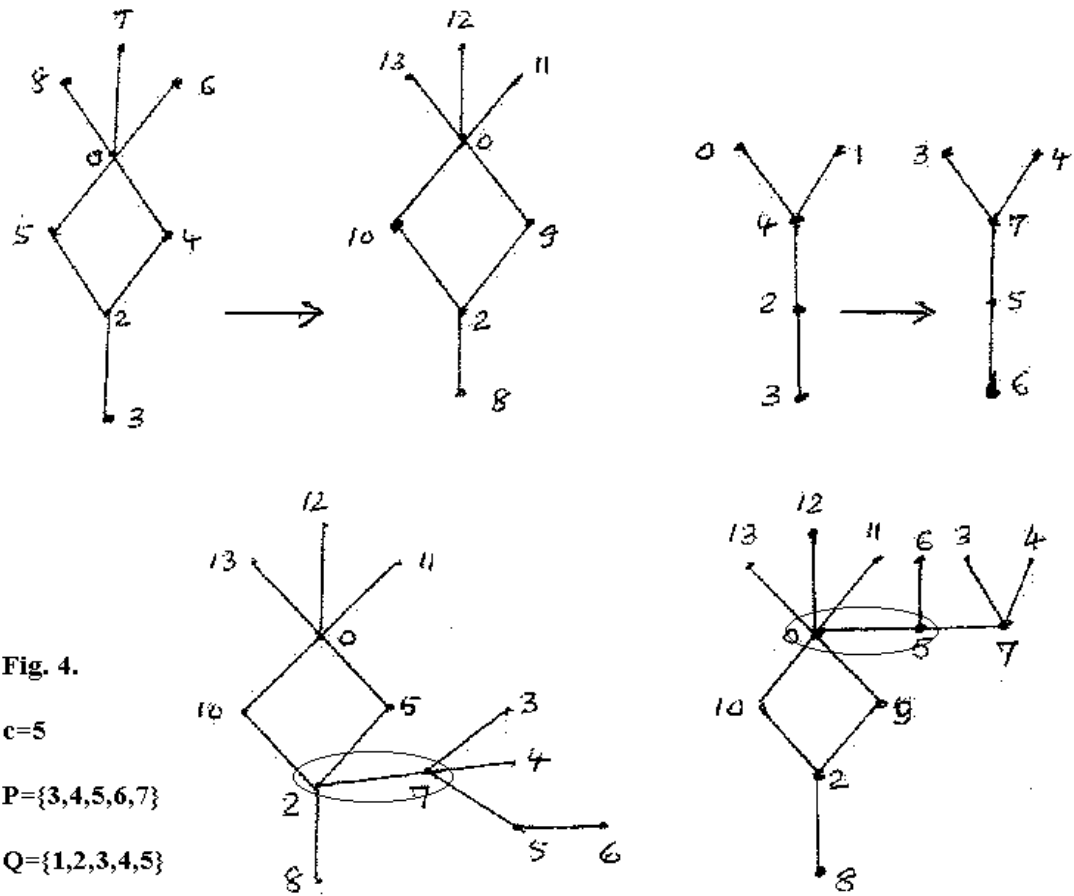

**Fig. 4.**

c=5

P={3,4,5,6,7}

Q={1,2,3,4,5}

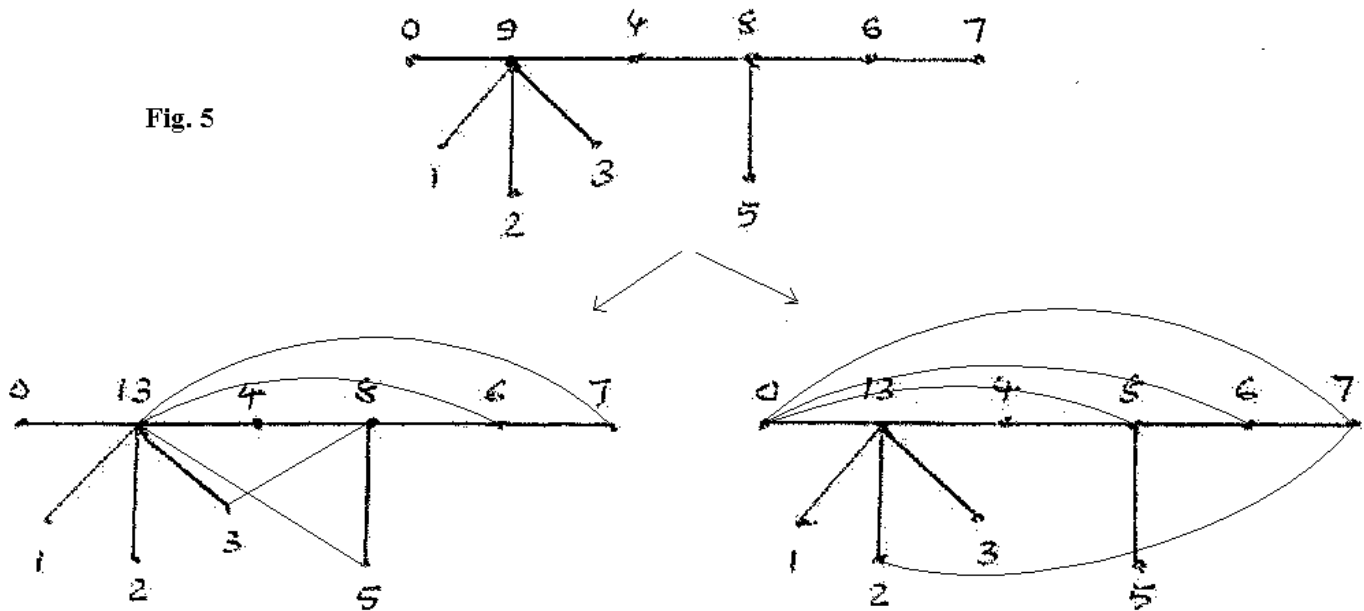

**Fig. 5**

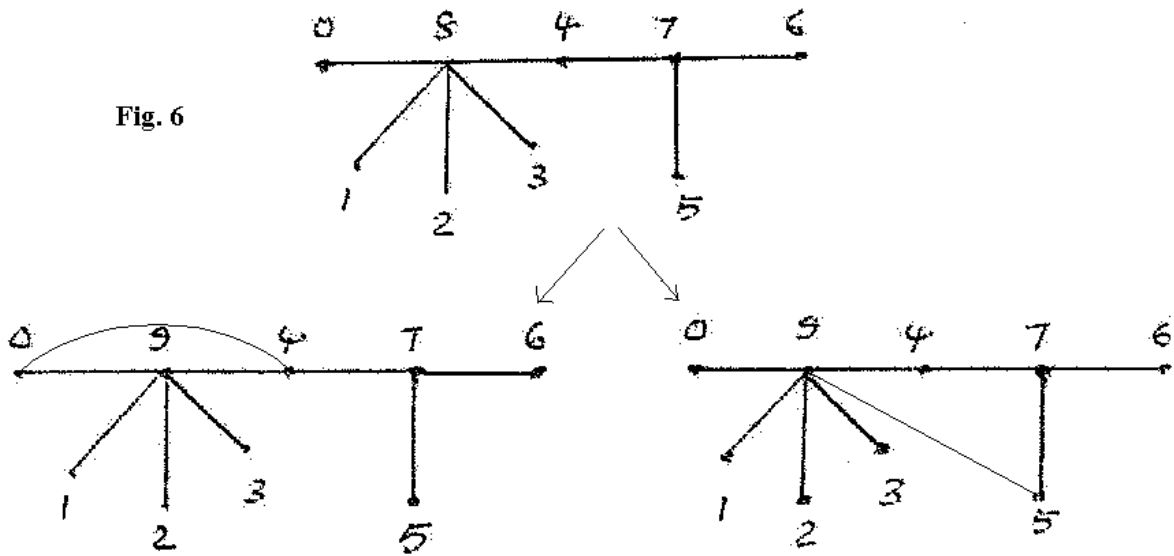

**Fig. 6**

**Possible Pairs: {(0,4), (1,5), (2,6), (3,7), (5,9)}**

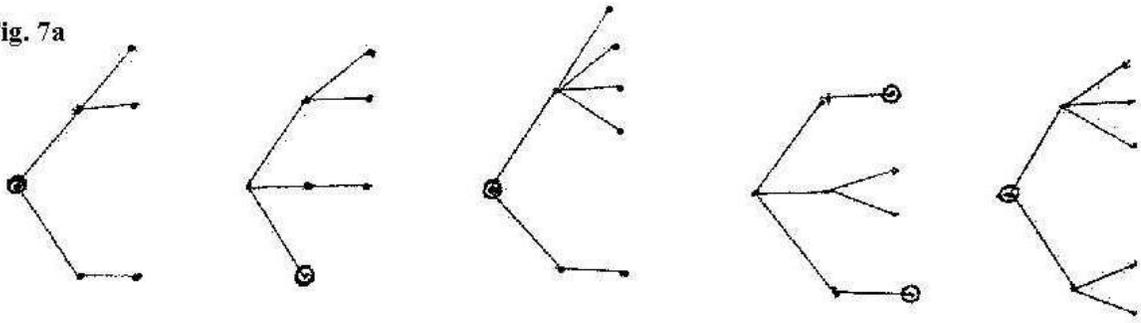

Fig. 7a

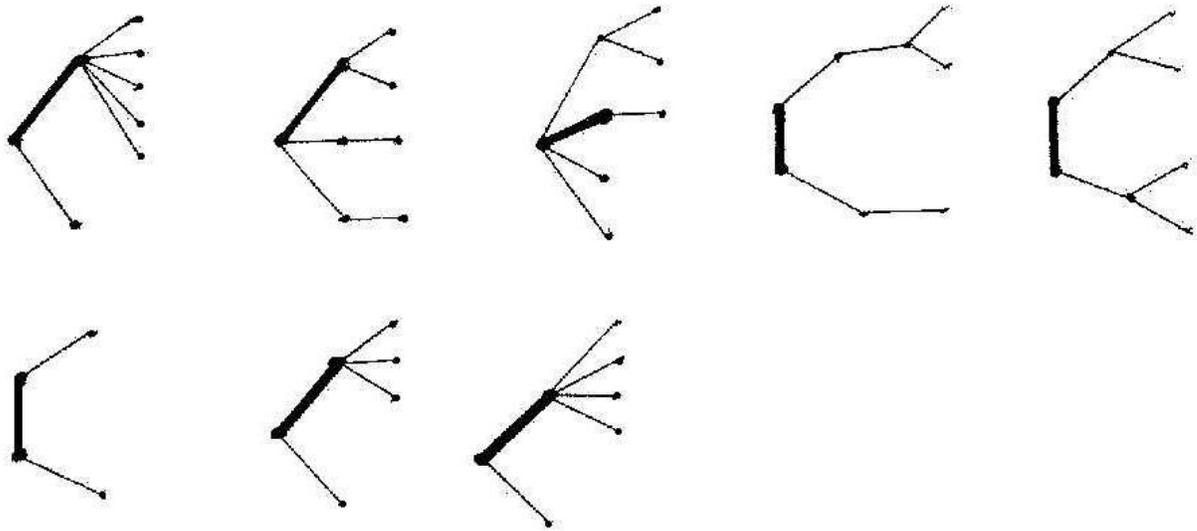

Fig. 7b

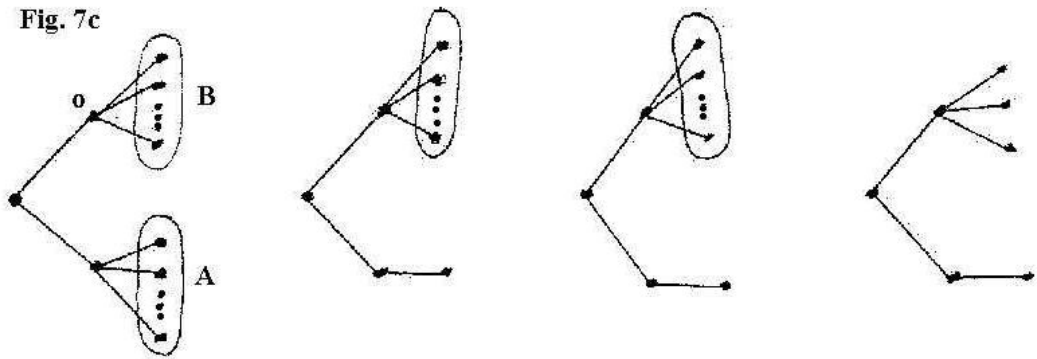

Fig. 7c

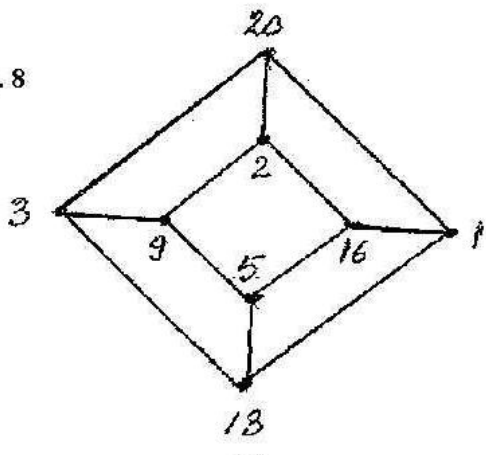

Fig. 8

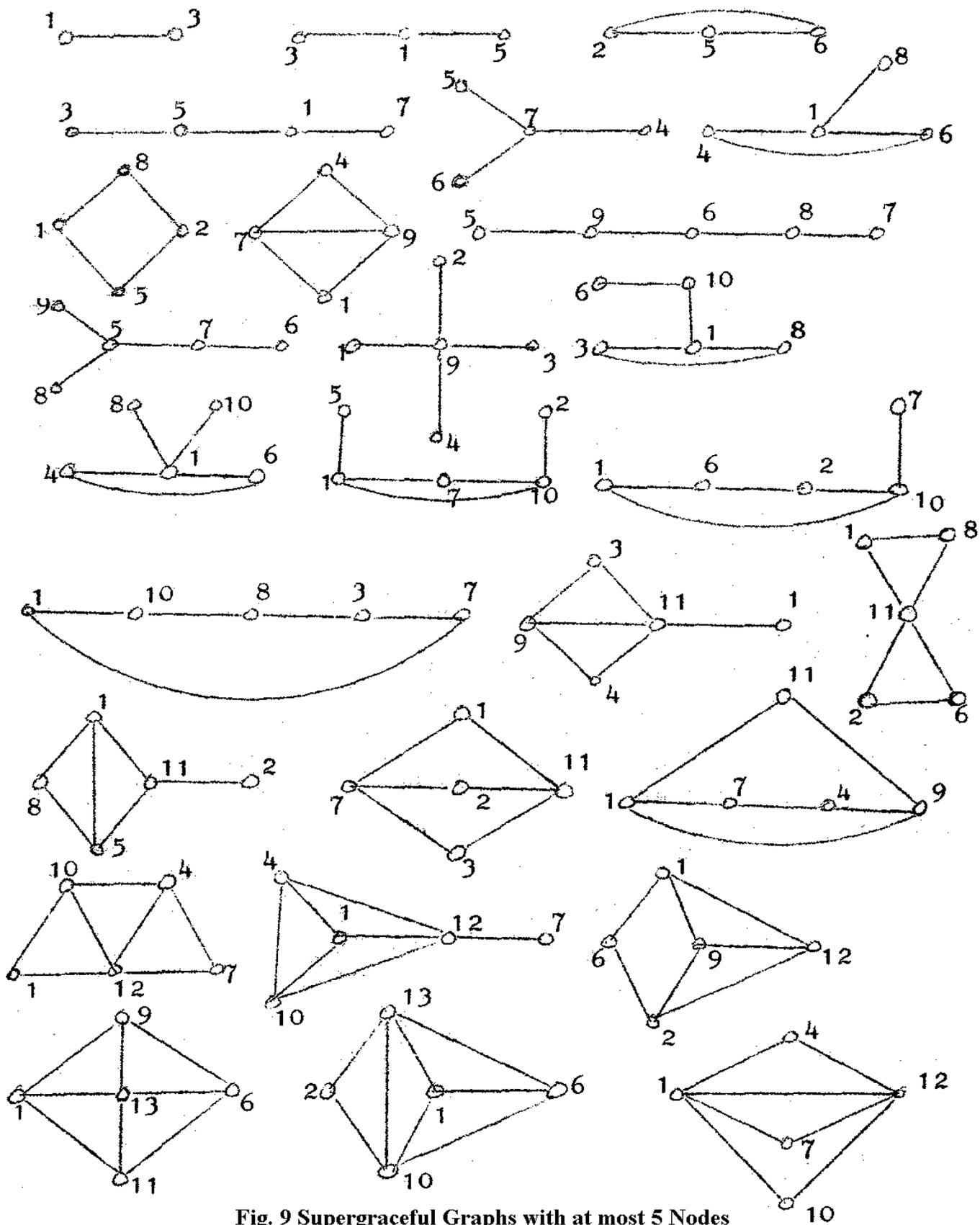

**Fig. 9 Supergraceful Graphs with at most 5 Nodes**

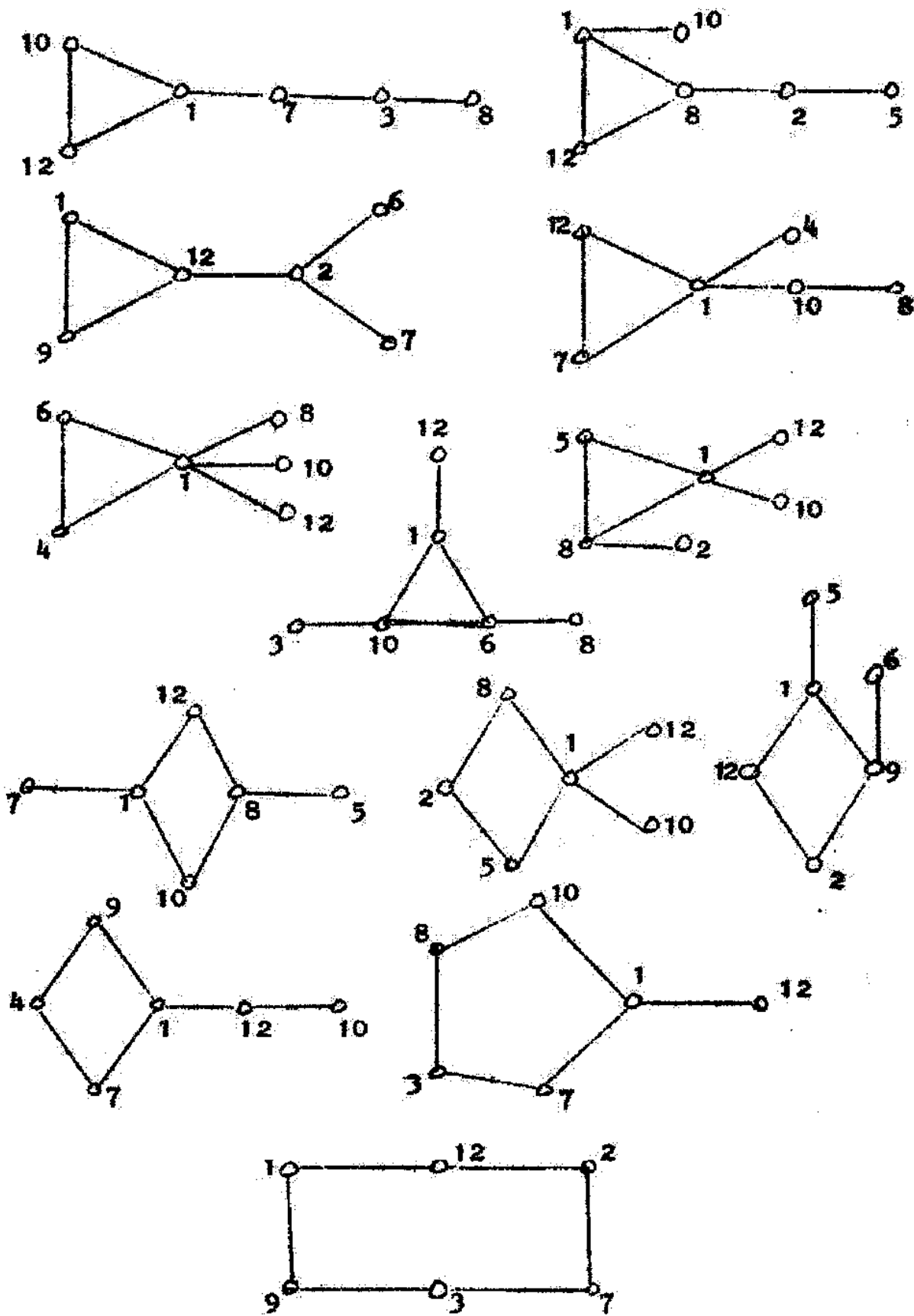

**Fig. 10. Unicyclic Graphs of Order 6 with a Total Numbering**

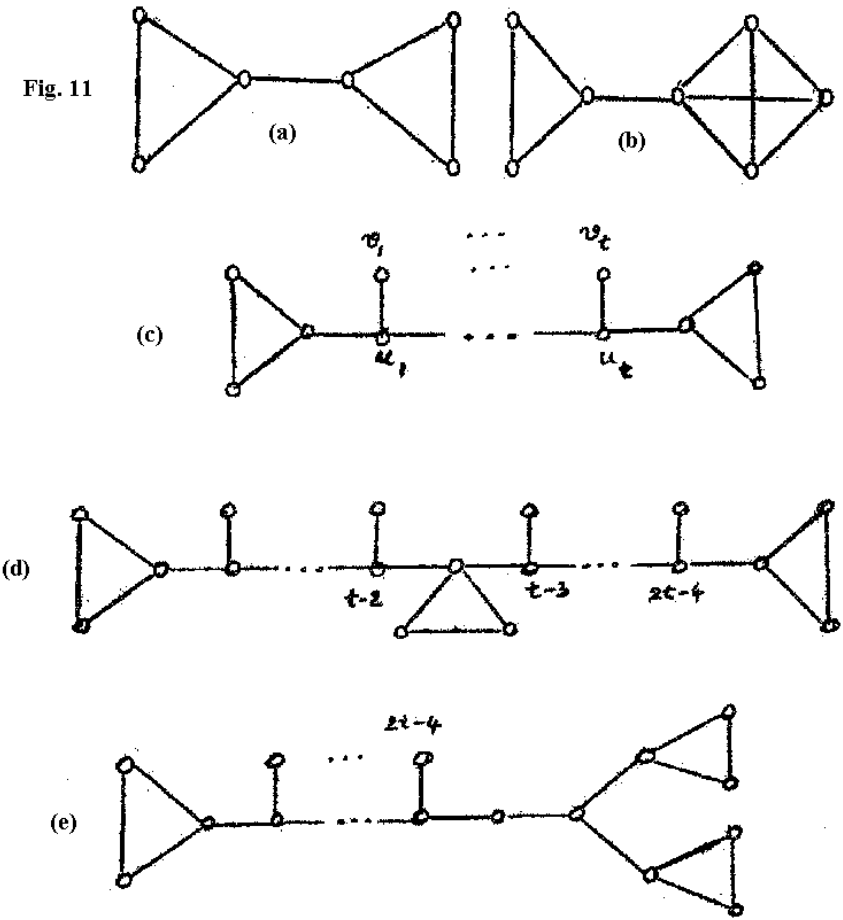

Fig. 11

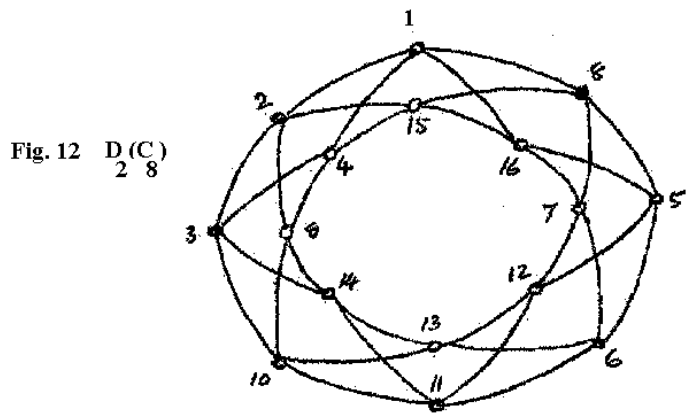

Fig. 12  $D_2(C_8)$

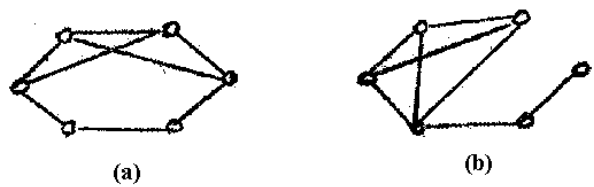

Fig. 13